\documentclass[review]{elsarticle}
\usepackage[utf8]{inputenc}
\usepackage{mathtools}
\usepackage[english]{babel}
\usepackage{amsthm}
\usepackage{amsmath}
\usepackage{amssymb}
\usepackage{comment}
\usepackage{url}
\usepackage{bbm}
\usepackage{bm}
\usepackage{subcaption}
\usepackage{multirow}
\usepackage{relsize}
\usepackage{color}
\usepackage{wrapfig}
\usepackage[export]{adjustbox}
\usepackage{tabulary}
\usepackage{booktabs}
\usepackage{color, soul}

\newcounter{stepno}

\newtheorem{lemma}{Lemma}[section]

\newtheorem{corollary}{Corollary}

\newcommand{\E}{\mathbb{E}}
\newcommand{\bn}{\begin{eqnarray}}
\newcommand{\en}{\end{eqnarray}}
\newcommand{\bns}{\begin{eqnarray*}}
\newcommand{\ens}{\end{eqnarray*}}

\newcommand{\defvarbegin}{\begin{quotation}\vspace{-15pt}\begin{tabbing}}
\newcommand{\defvarend}  {\end{tabbing}\vspace{-10pt}\end{quotation}}

\newcommand{\bnarray}{\begin{equation}\begin{array}{rcll}}
\newcommand{\enarray}{\end{array}\end{equation}}

\newcommand{\barr}{\begin{array}}
\newcommand{\earr}{\end{array}}

\newcounter{cnum}

\newcommand{\beginalg}{\setcounter{stepno}{1}
                \begin{list}{\bf Step~\arabic{stepno}}
                         {\usecounter{stepno}\settowidth{\labelwidth}{\bf Step~9m}
                \addtolength{\leftmargin}{2\parindent}}
                }
\newcommand{\eg}{\end{list}}

\def \define{\begin{quote}\begin{itemize}}
\def \enddefine{\end{itemize}\end{quote}}

\newlength{\boxedparwidth} 
  {\begin{center} \begin{tabular}{|@{\hspace{.15in}}c@{\hspace{.15in}}|}
                 \hline \\ \begin{minipage}[t]{\boxedparwidth}
                 \setlength{\parindent}{0.0in}}%
   {\end{minipage} \\ \\ \hline \end{tabular} \end{center}}
\newcounter{example}
\newcommand{\argmax}{{\rm arg}\max}

\def \xtilde{{\tilde x}}

\def \Btilde{{\tilde B}}

\def \Etilde{{\tilde E}}

\def \Rtilde{{\tilde R}}
\def \Stilde{{\tilde S}}

\def \ahat{\hat a}

\def \Lhat{\hat d}

\def \fhat{\hat f}

\def \khat{\hat k}

\def \yhat{\hat y}

\def \Ahat{\hat A}

\def \Dhat{\hat D}

\def \Hhat{\hat H}

\def \Lhat{\hat L}
\def \Khat{\hat K}

\def \Qhat{\hat Q}

\def \Uhat{\hat U}

\def \Pbb{\mathbb P}

\def \phihat{\hat \phi}

\def \hbar{\bar h}

\def \pbar{\bar p}

\def \Ncal{{\cal N}}

\def \Wcal{{\cal W}}
\def \Xcal{{\cal X}}
\def \Ycal{{\cal Y}}
\def \Zcal{{\cal Z}}

\def \xbf{{\mathbf{x}}}
\def \fbf{{\mathbf{f}}}
\def \cbf{{\mathbf{c}}}
\def \Mbf{{\mathbf{M}}}
\def \Vbf{{\mathbf{V}}}
\def \Abf{{\mathbf{A}}}

\biboptions{authoryear}

\graphicspath{ {./images/} }
\begin{document}

\begin{frontmatter}

\title{An Information-Collecting Drone Management Problem for Wildfire Mitigation}

\author[EEPU]{Lawrence Thul\corref{mycorrespondingauthor}}
\cortext[mycorrespondingauthor]{Corresponding author}
\ead{lathul@princeton.edu}
\tnotetext[t1]{This work is supported by the Air Force Office of Scientific Research (Award Number FA9550-19-1-0203).}
\author[ORFE]{Warren Powell}

\address[EEPU]{Department of Electrical Engineering, Princeton University, Princeton, New Jersey,USA}
\address[ORFE]{Department of Operations Research and Financial Engineering, Princeton University, Princeton, New Jersey,USA}

\begin{abstract}
We present a formal mathematical multi-agent modeling framework for autonomously combating a wildland fire with unmanned aerial vehicles.  The problem is formulated as a collaboration between a drone and a helicopter equipped with a tanker.  The modeling solutions are designed to capture the communication between agents and the information processes between the agents and their environment.  The drone is used to make partial observations and use them to update probability distributions over the uncertain state of the world.  We design a parameterized direct lookahead approximation policy to guide the drone through the region and promote active learning.  The helicopter is used to extinguish the fire.  We design a helicopter policy which anticipates the zones which will burn through the belief modeling and extinguish them.  We simulate the modeling and policy solutions using a wildland fire simulation in an area of Northern California.
\end{abstract}


\end{frontmatter}

\section{Introduction}
\label{sec:introduction}
Wildland fires are a necessary component of the Earth's complex climate; however, under certain conditions they can be very costly to the economy, natural habitats, and human lives.  As the climate changes, the dry seasons are longer, the fires grow bigger, and resources deplete faster.  Wildland fires have become more prevalent in remote areas, such as the arctic, where it is difficult and dangerous to send resources.  The proper management of resources can lead to better active planning and firefighting techniques.  Multi-agent autonomous systems can add value to this problem because they can quickly dispatch into areas to contain fires which mitigates the risks to human lives and reduces operating costs.  In this paper, we present a multi-agent framework for collaborative firefighting between an unmanned aerial vehicle (drone) and a helicopter fitted with a tanker to contain the fire under a critical time horizon.

The different types of agents will work together to achieve the joint goal of mitigating the cost of burning fuel while also performing individual tasks.  There are a number of complexities involved with constructing an autonomous system for wildland fire fighting.  For example, the system must maintain a belief about the locations of the fire, manage battery life, and decide the best areas to extinguish.
There are several methods of observing the existence of a fire such as observation towers, satellite imagery, or smoke on the horizon; however, these methods only provide high-level information about the approximate location of the fire.  It would be much more insightful to obtain localized intensity readings to guide the extinguishing decisions.  Drones provide a method to obtain high precision readings about the local information about fire intensity in real time.  However, drones can only be in one place at a time so the observations will only provide partial information about the state of the region.  The drone communicates the updated belief about the state of the fire with the helicopter.  The helicopter is used to directly extinguish the fire at each decision epoch.  The stochastic dynamic processes describing the fire spread are highly uncertain and can cover large spatial areas fairly quickly, so it is imperative that the helicopter anticipate which decisions will make the largest impact on mitigation.  Additionally, the drone has a limited battery life, so it must be managed to ensure it returns home.  The modeling and planning algorithms must incorporate real-time logic to update probability distributions about the state of the fire, forecast into the future, and actively plan routes through the region.

In this work, we extend the unified modeling framework for stochastic optimization described in \cite{PowellEJOR} to handle multiple agents.  This framework is unique because it emphasizes the notion of modeling the problem first, then performing a search over the four classes of policies to design a solution.  There are various characteristics which make this problem difficult.  Firstly, the state of the world is not fully observable, so there must be probability distributions over the uncertain state variables (belief states).  The spatiotemporal stochastic dynamic processes are highly uncertain which makes predicting when and where a new zone will ignite can be challenging.  The drone must manage the value of observations with the physical constraints of its position, speed, and battery life.  Additionally, the size of the region is large and making effective decisions and predictions requires modeling the fire at a high spatial resolution.  Hence, the solution strategies will encounter three curses of dimensionality which will limit the set of feasible solutions.  The state space, action space, and observation spaces are very high dimensional so optimal policies are infeasible to find; however, there is structure in this problem we exploit to find effective strategies to mitigate the spread.

This paper makes the following contribution:
\begin{itemize}
    \item We present the first formal mathematical model for a multi-agent wildfire model using the unified framework,
     \item We design a hybrid CFA-DLA policy for the helicopter which trades off between exploiting where to extinguish fire now versus the cost of zones where it may spread to,
    \item We design policies for the drone agent which make observations through active learning while managing a physical state using direct lookahead strategies.  We present a Thompson sampling DLA and a parameterized interval estimation DLA and compare the results of each,
    \item We demonstrate that the interval estimation direct lookahead drone policy in conjunction with the hybrid CFA-DLA helicopter policy performs best when we apply the multi-agent framework in a simulator designed using data from a small region in Northern California with frequent fires.
\end{itemize}

This paper is organized as follows.  Section 2 reviews the literature about wildland fire planning and modeling, and then covers the literature about stochastic optimization frameworks in routing and transportation problems.  Section 3 provides a detailed description of the problem setting. Section 4 presents the formal mathematical model for the multi-agent system.  The section includes the environment model, the drone model, and the helicopter model.  Section 5 designs policies for making decisions, including the information collection policy for the drone and the implementation policy of the helicopter.  The first part of the section presents the heuristic helicopter policy, and the second part presents the policies for making learning decisions with the drone.  Performance results for the policies using our wildland fire simulator are presented in Section 5.  Section 6 draws conclusions and discusses future work.

\section{Literature Review}
\label{sec:lit_review}
This section is decomposed into two general parts.  First, we outline the relevant wildland fire literature.  Then, the relevant stochastic planning literature with regards to vehicle routing and resource allocation under uncertainty. \cite{jain2020review} reviews artificial intelligence and data science methods related to wildland fire science and management.  There are many different topics under the domain of wildland fire management, but we narrow our review to the relevant topics in spatiotemporal wildland fire prediction, monitoring wildland fires, and resource allocation and planning.

There are several methods for wildland fire prediction.  Some focus on physics based models, while others rely on simulation methods.   \cite{sullivan2009wildlanda} outlines the physics based models used in the literature. \cite{andrews2018rothermel} provides a comprehensive review of Rothermel’s fire surface spread model, which incorporates dynamic fuel models, weather, and topography to predict the rate of spread. There are also methods which take a mathematical modeling approach by using computer simulations.  \cite{lopes2002firestation} describes a software simulation system for forecasting wildfires across complex topographies.  \cite{finney1998farsite} provides details of the FARSITE simulator which has been used in various scientific studies.

Many statistical and machine learning techniques have been used to develop forecasting methods for wildfires.  \cite{ zhou2019ensemble} and \cite{srivas2016wildfire} use ensemble Kalman filtering to make fire predictions.  \cite{zheng2017forest} uses cellular automatons to learn a model for wildland fire prediction.  \cite{subramanian2017learning} uses reinforcement learning to predict fire growth from satellite data.  A common approach in the recent literature is to use deep neural networks to make predictions (e.g. \cite{lattimer2020using}, \cite{radke2019firecast}, \cite{vakalis2004gis}).

During a real-time wildland fire event, it is imperative to understand the current state of the fire.  There are various methods for observing fires such as fire watch towers, satellites and unmanned aerial vehicles (UAVs).  \cite{bao2015optimizing} solves a multi-objective genetic algorithm for optimizing the location of watch towers.  \cite{markuzon2009data} generates predictions from satellite data.  \cite{merino2012unmanned} and \cite{zhao2018saliency} use UAVs with on-board camera imaging to monitor fires in real-time systems via thermal image feature matching and deep learning techniques, respectively.   \cite{beard2006decentralized} and \cite{casbeer2005forest} use multiple UAVs to collaboratively collect the most valuable information to monitor a wildfire.  \cite{julian2019distributed} uses a deep reinforcement learning algorithm to control the flight path of a UAV to monitor a wildfire in a partially observable setting.

We also outline the relevant literature related to wildland fire resource allocation and planning.  There are preemptive measures that could be taken to prevent wildfires and there are active response strategies which must be deployed to combat them.  We will focus on the latter in this review.  See \cite{thompson2011uncertainty} for a comprehensive review of uncertainty and risk management during wildland fire events.  \cite{calkin2011real} provides a system for real-time risk assessment for supporting resource allocation decisions.  \cite{mcgregor2016fast} uses an MDP model and MCTS algorithm for simulating fire suppression policies.  \cite{bertsimas2017comparison} provides a rolling horizon Monte Carlo Tree Search (MCTS) algorithm and applies it to solving a tactical wildfire resource allocation problem.  \cite{dimopoulou2004towards} provides an integrated framework for the control of resources during a forest fire which uses simulation and optimization based solutions.  The onset of a wildfire is the most critical time to mitigate the spread before it becomes impossible to contain.  \cite{maclellan1996basing} optimizes the placement of air-tankers across firefighting bases in Ontario Canada to ensure they are properly distributed throughout the region.  \cite{hodgson1978location} optimizes a single extinguishing route of an air-tanker with an integer programming solution.  \cite{podur2007simulation} provides simulations to allocate resources after a fire has grown very large.
Problems related to stochastic dynamic vehicle routing for autonomous systems arise in various areas of the sequential decision modeling literature.   \cite{ritzinger2016survey} provides a survey on solution strategies for dynamic vehicle routing problems.  The problem becomes more challenging when there is uncertainty about the state of the world.  Vehicle routing problems when the underlying problem is a partially observable Markov decision process (POMDP) have been solved via stochastic control solutions (e.g. \cite{liu2019vehicle}).  Decision-making under our partially observable problem setting can be classified as a partially observable Markov decision process (POMDP).  The POMDP approach is a common strategy for modeling problems with uncertainty around the true underlying state variable of the system, but they have computational limitations because it is PSPACE-complete (e.g \cite{pineau2006}).  Some approximation strategies have been proposed in the literature.  For example, \cite{ross2008} derives online planning algorithms for the POMDP problem, \cite{pineau2006} derives point-based solutions for solving large POMDPs,  and \cite{roy2005} reduces the dimensionality of the belief state by projecting it into a lower dimensional space.  In our problem setting, none of these strategies will be able to handle the high dimensionality of the state space, action space, and observation space.

Any solution to a dynamic vehicle routing problem under uncertainty can be categorized as a stochastic optimization problem.  There are many different communities for solving sequential decision and stochastic optimization problems; however, the policies used for making decisions can be aggregated into four classes of policies outlined in \cite{PowellEJOR}.  We formulate the solutions in this paper by drawing our approximate problem strategies from the four classes of policies.  We have discussed papers in the literature which have considered the problem of monitoring wildfires with single and multiple agents in a partially observable setting.  We have also discussed papers which have considered the problem of allocating resources to extinguish fire in an active event.  This is the first paper to consider joint firefighting and fire monitoring agents in a partially observable setting.   We present the first multi-agent mathematical formulation of a wildland fire problem with the unified modeling framework.  We operate a helicopter for mitigating the spread of wildfire in conjunction with a drone for monitoring wildfire in a partially observable setting.  We present a parameterized multi-stage rolling horizon optimization for learning with the drone which can handle the physical constraints of the problem.

\section{Problem Description}
\label{sec:problem_desc}
In this work, we consider the joint collaboration between an unmanned aerial vehicle (drone) and a helicopter equipped with a tank to mitigate the spread of a wildland fire before the end of time horizon, $T$.  We assume the wildland fire has been spotted in a region split into a set of zones, $\Zcal$.  Each zone, $z\in \Zcal$, has a fixed amount of fuel particulates, $\eta_z \in \mathbb{N}$.  At each time $t$, the fuel in the zone can be partitioned into a portion of healthy fuel $H_{tz}$, a portion of burning fuel $Q_{tz}$, and a portion of dead fuel $D_{tz}$.  The rate the fuel burns is determined by a set differential equations parameterized by static fuel characteristics, $\Xi_z$, which are known to the agents (See Appendix).  Additionally, the rate the fuel spreads between zones is determined by the fuel characteristics, the wind speed $U_t$, the wind direction $\phi_t$, and the topography of the map, $\kappa_z$.  There is also a cost, $r_z$, associated with the fuel in each zone which corresponds to the value of houses, roads, or endangered habitats versus less valuable areas such as grasslands or shrubs.

We assume there is a home zone $z_0\in \Zcal$ which the helicopter and drone are dispatched from at the beginning of the event.  The drone is capable of observing fire within a field of view $\rho^{obs}$ from its physical location, $R_{t}^{drone} \in \{0,1\}^{|\Zcal|}$, where
\bns
R_{tz}^{drone} = \begin{cases}
    1 & \text{if the drone is in zone $z$ at time $t$},\\
    0 &  \text{otherwise},
\end{cases}
\ens
but it has a limited amount of battery life, $R_t^{batt}\in \mathbb{R}$.  There is enough battery life to last the time horizon $T$, but the drone must return home at the end of the horizon to make sure it is not lost.  We assume there is a maximum speed the drone can travel and the drone consumes a constant amount of energy per unit time.  At each arrival to a destination node, the drone must stop to capture and process the new intensity readings.  With the intensity readings, the drone can estimate which zones are burning and how much fuel in that zone is ignited with relatively low noise.

\begin{figure}[ht]
\includegraphics[width=0.75\textwidth, center]{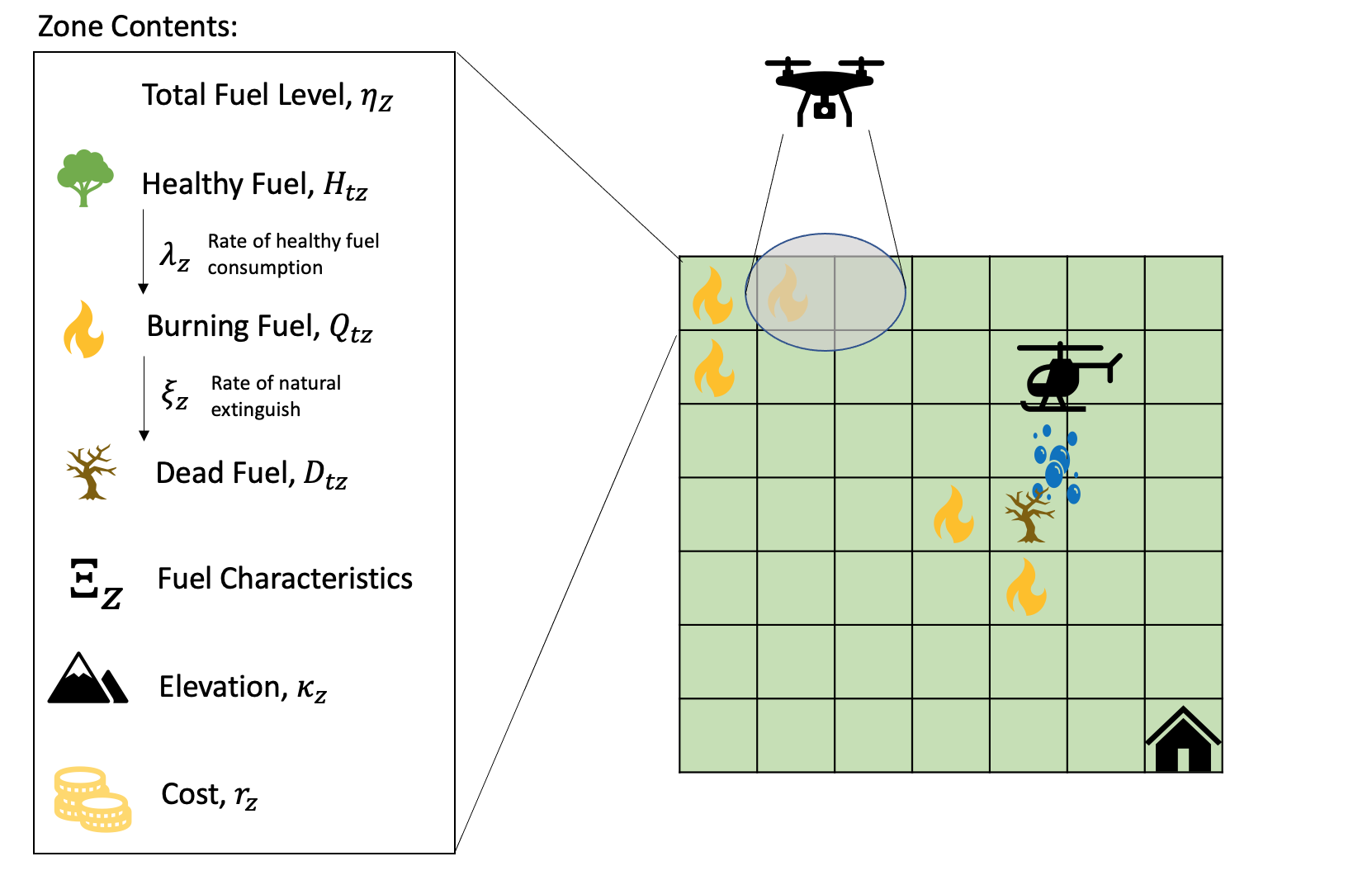}
\caption{The box on the left gives descriptions of the features of a zone $z\in\Zcal$ at time $t$.  The illustration depicts a snapshot of the problem setting. }
\label{fig: illustration}
\end{figure}

The helicopter is responsible for directly extinguishing the fire by releasing fire retardant from its tank.  It is common for helicopter tankers to refill the fire retardant from a nearby water-source, which we assume can be performed in under ten minutes. The helicopter has a larger fuel capacity than the drone, so its fuel level is negligible over the time horizon.  For both agents, we assume the decision epochs are always every 10 minutes.  The helicopter is capable of refilling under those time constraints and the drone must move and capture/process images.  However, the drone moves much slower than the helicopter so the distance it can travel in one time step is constrained.  Therefore, the goal of the problem is to design a decision-making strategy for each agent which can operate in real time and contain the fire under the constraints laid out in the problem description.  Figure \ref{fig: illustration} outlines the contents of each zone and presents an illustration to reinforce the general layout of the problem.

\section{Modeling Framework}
\label{sec:model}
In this section, we provide the multi-agent extension of the unified modeling framework for the problem described in section \ref{sec:problem_desc}.  There are three different perspectives to model: the environment agent, the learning agent (drone), and the implementation agent (helicopter).  The environment agent represents a dynamic state of the world which can be impacted by exogenous factors or other agents, but it is passive (makes no decisions).  There are also controlling agents which make decisions sequentially through time to learn about the environment or directly change it.  The set of collaborative agents will have a joint goal, but also have individual goals.

In this problem class, the joint goal is to mitigate the total cost of burning fuel, and each agent has an individual job to accomplish that joint goal.  The utility of the learning agent is to collect the most valuable observations.  The utility of the helicopter agent is to directly extinguish the fire.  The information process from the perspective of each agent will be different.  Throughout this paper we introduce notation to keep track of the information known to each agent.  Variables superscripted by $env$, $heli$, and $drone$ represent variables within the environment model, helicopter model, and drone model, respectively.  Variables with no accent are the true values and variables with bars are estimates of those values.  State variables superscripted by $x$ are post decision state variables from the helicopters extinguishing decision.  The expectation using the knowledge held within the belief state, $B_t$, through time $t$ are reflected by $\E[\cdot|B_t] = \E_t[\cdot]$.

The environment agent will be impacted by exogenous information (e.g. wind changes) and the helicopter decisions, $\xbf_t^{heli}$.  The drone will make new observations, update the belief state, then make a decision, $\xbf_t^{drone}$, to make a new set of observations at a new location.  Then, the helicopter uses the updated belief state to make a new extinguishing decision.  The following sequences show the information processes for each agent as they unfold:

\begin{itemize}
    \item Environment Agent;
    \bns
    (..., S_{t}^{env}, x_{t}^{heli}, S_{t}^{env,x}, W_{t+1}^{env},  S_{t+1}^{env}, x_{t+1}^{heli}, S_{t+1}^{env,x}, W_{t+2}^{env}, S_{t+2}^{env}, ...)
    \ens
    \item Learning Agent (Drone);
    \bns
    (..., S_{t}^{drone}, x_t^{heli}, S_{t}^{drone,x}, x_t^{drone}, W_{t+1}^{drone}, S_{t+1}^{drone}, x_{t+1}^{heli},  S_{t+1}^{drone,x}, x_{t+1}^{drone}, W_{t+2}^{drone}, S_{t+2}^{drone}, ...),
    \ens
    \item Implementation Agent (Helicopter);
    \bns
    (..., S_{t}^{heli}, x_t^{heli}, W_{t+1}^{heli}, S_{t+1}^{heli}, x_{t+1}^{heli}, W_{t+2}^{heli}, S_{t+2}^{heli}, ...).
    \ens
\end{itemize}


The following sections will model each agent using the unified framework.  The unified framework can model any sequential decision problem with five components: the state variable, the decision variable, the exogenous information, the transition function, and the objective function.

\subsection{Environment Agent}
\label{sec:environment}
The environment agent represents the dynamically changing environment.  We model the environment agent to have perfect knowledge about the state of the world.  Since a passive agent does not make decisions, then it will have no decision variable or objective function.

Recall, we assume the zone has a fixed amount of fuel, $\eta_z$, and this fuel will remain in the region (in various states) throughout the time horizon $T$.  At time $t$, the fuel in each zone has a state consisting of:
\begin{itemize}
    \item $H_{tz} =$ the amount of healthy fuel at time $t$ in zone $z$,
    \item $Q_{tz} =$ the amount of burning fuel at time $t$ in zone $z$,
    \item $D_{tz} =$ the amount of dead or extinguished fuel at time $t$ in zone $z$.
\end{itemize}
The intensity of the fire in the zone, $Y_{tz}$ is proportional to the amount of fuel burning at $t$; hence, $Y_{tz} \propto Q_{tz}$.  Since the amount of fuel in the zone is fixed, then the fuel states of each zone exhibit the following property,
\bn
\label{eqn:fuel_prop}
H_{tz} + Q_{tz} + D_{tz} = \eta_z.
\en
Furthermore, it will be convenient to introduce an indicator variable which represents whether a zone has ignited given by,
\bn
\label{eqn:Kdef}
K_{tz} = \begin{cases}
      1 & Y_{tz}>0, \\
      0 & otherwise.
   \end{cases}
\en

At each time $t$, we assume a helicopter is in the region and sequentially dropping fire retardant on a set of zones.  The helicopter makes a decision, $\xbf_{t}^{heli} \in \{0,1\}^{|\Zcal|}$, to drop fire retardant at a specified zone, where
\bns
x_{tz}^{heli} = \begin{cases}
      1 & \text{extinguish zone $z$}, \\
      0 & otherwise,
   \end{cases}
\ens
constrained by $\sum_{z\in\Zcal} x_{tz}^{heli} = 1$.

When the fire retardant is released on zone $z$, it extinguishes the fire completely in that zone as well as some neighboring zones $z' \in \Zcal_{t}^{heli}$, which is given by,
\bns
\Zcal_{t}^{heli} = \{ z' \in \Zcal: ||z' - x_{tz}^{heli}|| \leq \rho^{heli}\},
\ens
where $\rho^{heli}$ is the extinguishing radius and $||\cdot||$ is an operator computing the distance between two zones.  Additionally, assume that once a zone is extinguished it will not reignite during the time horizon $T$.

The behavior of fire as it spreads throughout the region is highly stochastic.  To capture the spatial behavior of the model, assume the set of zones in the region are nodes in a directed graph.  Between time $t$ and $t+1$, assume that an edge $(z',z)\in \Zcal \times \Zcal$ is added to the graph if a spreading fire traverses the distance $\|z-z'\|$ from zone $z'$ and ignites $z$.  The random distance the fire will spread from $z'$ to $z$ is given by the random variable $\Lhat_{t+1,zz'}$ and it depends on multiple factors such as the heterogeneous fuel in the zones, the wind speed and direction, and the slope of the land.  Therefore, the adjacency matrix for the graph at time $t+1$, $\mathbf{\Ahat}_{t+1}$ is a random binary matrix.  An entry of the adjacency matrix, $\ahat_{t+1,zz'}$ is $1$ if an edge connects zone $z'$ to $z$, and $0$ otherwise.  If $z \in \Zcal_{t}^{ext}$, then $\ahat_{t+1,zz'} = 0$  $\forall z' \in \Zcal$.  In other words, it is removed from the graph because there are no directed edges in or out of it.

The entries of the adjacency matrix can be used to describe the random dynamics for fire spreading between zones.  The updating equations for the fire spreading dynamics are given by,
\bn
\label{eqn:env_ext_update}
\Zcal_{t}^{ext} &=& \Zcal_{t-1}^{ext} \cup \Zcal_{t}^{heli}(x_t^{heli}), \\
\label{eqn:env_Kx_update}
K_{tz}^{x} &=& K_{tz}  \mathbf{1}_{[\Zcal \setminus \Zcal_{t}^{ext}]}, \\
\label{eqn:env_Kt1_update}
K_{t+1,z} &=& 1 - \prod_{z' \in \Zcal} \left(1 - \ahat_{t+1,zz'} K_{tz'}^x\right),
\en
where equation \eqref{eqn:env_ext_update} adds the newly extinguished zones to the set of extinguished zones, equation \eqref{eqn:env_Kx_update} sets all indicator variables corresponding to the extinguished zones to zero, and equation \eqref{eqn:env_Kt1_update} updates the value of the indicator variables at $t+1$ by passing the values through the binary adjacency matrix.  Furthermore, there exists a set of differential equations describing the fuel states as they transition through time.  We leave them as a general function given by,
\bn
\label{eqn:transition_general}
(H_{t+1,z}, Q_{t+1,z}, D_{t+1,z}) &=& f^{env}(S_t^{env}, \xbf_t^{heli}, W_{t+1}^{env}),
\en
where $S_t^{env}$ is the state variable of the environment agent and $W_{t+1}^{env}$ is the exogenous information for the environment agent.  The equations are known to the environment agent, but the explicit forms of these equations are only needed for the simulator defined in Appendix \ref{app:Rothermel}.

The environment agent model can be modeled using the five components of the unified framework.  Since the environment agent is a passive system it does not make any decisions; hence, it has no objective function.  However, the decisions made by the controlling agents have an impact on the physical system.

\paragraph{Environment State Variable} The state variable for the environment agent describes the dynamically changing variables needed to compute the transition functions.  The pre-decision environment state variable is given by,
\bns
S_{t}^{env} = \left(\left(Q_{tz}, H_{tz}, D_{tz}, K_{tz} \right)_{z\in\Zcal}, \Zcal_{t-1}^{ext}\right).
\ens

After each decision is made, the post-decision state variable describes the environment state after the controlling agent's have made decisions, but there is not yet any new information available.  In other words, the decisions have been implemented, but it is still time $t$.  The post-decision environment state is given by,
\bns
S_{t}^{env,x} = \left(\left(Q_{tz}^x, H_{tz}^x, D_{tz}^x, K_{tz}^x\right)_{z\in\Zcal}, \Zcal_{t}^{ext}\right).
\ens

The initial state variable consists of the initial state of each dynamic state variable, as well as the static variables through time which are kept in the system as latent variables.  The initial state variable is given by,
\bns
S_{0}^{env} = \left(Q_{0z}, H_{0z}, D_{0z}, K_{0z}, \eta_z,  \Xi_z\right)_{z\in\Zcal}
\ens
where,
\bns
\Xi_z &=& \text{Set of static fuel parameters in $z$}, \\
\eta_z &=& \text{the total fuel level for zone $z$}.
\ens
The initial state consists of the initial dynamic state variables plus static parameters which are not changing over the time horizon.  We omitted the extinguished set since it is empty at time $0$.

\paragraph{Environment Exogenous Information}
The exogenous information describes the information streaming into the environment model between time $t$ and $t+1$ from exogenous sources or decisions made by other agents in the system.  The information which arrives between time $t$ and $t+1$ is given by,
\bns
W_{t+1}^{env} = \left(\mathbf{\Ahat}_{t+1}, \Zcal_t^{heli} \right)
\ens
where $\mathbf{\Ahat}_{t+1}$ is the random binary adjacency matrix for the graph.  The random edges within this matrix are realizations of random distances produced by factors such as wind speed, wind direction, fuel characteristics and terrain.  $\Zcal_t^{heli}$ is produced from the helicopter decision and removes the extinguished nodes from the graph.

\paragraph{Environment Transition Function}
The transition function consists of the set of equations for updating each state variable at $t+1$.  For this problem, these equations can be split into two sets of equations.  There are the set of equations for updating the fuel state levels within zones and the set of equations for the fire surface spreading between zones.  The transition equations for the interzone spreading are given by equations (\ref{eqn:env_ext_update}) - (\ref{eqn:env_Kt1_update}) and the transition equations for the intrazone fuel dynamics are given by equation \eqref{eqn:transition_general}.

\subsection{Drone Agent}
\label{sec:drone_framework}
The drone is the learning agent and will make decisions to observe subsets of the region to inform the helicopter model where to extinguish.  When the learning agent makes a new observation, the probability distribution over the unknown state of the world will be updated to reflect the new knowledge.  In the following sections, we introduce the belief model and the updating equations for it.  Then, the canonical model is constructed for the drone with the five components of the unified framework.

\subsubsection{Belief Model}
\label{sec:belief_model}
The belief model represents a probability distribution over the parameters of the environment model which the controller cannot observe perfectly at each time step.  To formulate the belief model, there are three important components to outline:
\begin{itemize}
    \item the initial belief assumptions,
    \item the belief state,
    \item the belief state update.
\end{itemize}
The initial belief state assumptions describe anything the controlling agent is asserting to be true when modeling the environment.  The belief state describes the parameters of the probability distribution over the environment state variables at time $t$.  The belief state update describes the set of equations used to update each of the belief state variables after new assumptions have been made.

\paragraph{Initial Belief Assumptions:} There are structural characteristics which the controllers assume about the environment agent to build probabilistic models.  It is assumed the controlling agents do not know the fuel states $Q_{tz}$, $H_{tz}$, or $D_{tz}$ for all zones $z\in \Zcal$, so they will have corresponding belief state variables.  Additionally, $K_{tz}$ is an indicator variable for whether the zone has ignited and is unknown for each zone but can be observed through the observations.  The zones which have been permanently extinguished through time $t$ are known perfectly to the controllers, and the set of fuel characteristics in the initial state are known perfectly to the controllers.  It is reasonable to make these assumptions because the fuel characteristics are mapped by ecologists and forest services (e.g. \url{https://www.landfire.gov/}).  The drones will assume a set of environment state transition functions, but the arguments are random variables so the forecasted states will be random variables.

The spread of the wildfire is modeled through equations \eqref{eqn:env_ext_update} -\eqref{eqn:env_Kt1_update}, and we assume the distances, $\Lhat_{t+1,z}$, are random variables with directionally dependant exponential distributions.  The fuel dynamics within a zone are assumed to be a set of difference equations.  The process as the fuel is consumed is fairly rapid in the beginning and tails off as the fire smolders.  The fuel transitions are assumed to follow the following equations,
\bn
\label{eqn:env_Qtx_update}
Q_{tz}^x &=& Q_{tz} K_{tz}^x, \\
\label{eqn:env_Htx_update}
H_{tz}^x &=& H_{tz}, \\
\label{eqn:env_Dtx_update}
D_{tz}^x &=& D_{tz} + Q_{tz} (1 - K_{tz}^x), \\
\label{eqn:env_Qt1_update}
Q_{t+1,z} &=& Q_{tz}^x - \xi_z Q_{tz}^x + \lambda_z Q_{tz}^x H_{tz} + Q^{init}_{z} (K_{t+1,z} - K^x_{tz}) , \\
\label{eqn:env_Ht1_update}
H_{t+1,z} &=& H_{tz}^x - \lambda_z Q_{tz}^x H_{tz} - Q^{init}_{z} (K_{t+1,z} - K^x_{tz}), \\
\label{eqn:env_Dt1_update}
D_{t+1, z} &=& D_{tz}^x + \xi_z Q^x_{tz} ,
\en
where $\lambda_z$ is the rate at which healthy fuel is ignited by the burning fuel, $\xi_z$ is the rate at which the burning fuel naturally extinguishes, and $Q^{init}_{z}$ is the initial amount of healthy fuel lit upon ignition.  Equations \eqref{eqn:env_Qtx_update} - \eqref{eqn:env_Dtx_update} provide the post-decision transitions and equations \eqref{eqn:env_Qt1_update} - \eqref{eqn:env_Dt1_update} provide the assumed transitions between time $t$ and $t+1$.

The drone has a field of view from a given position.  The set of observable zones at $t$ within its field of view is given by,
\bns
\Zcal_{t}^{obs} = \{z\in \Zcal: ||z - R_{tz}^{drone}|| \leq \rho^{obs}\},
\ens
where $\rho^{obs}$ is the distance in each direction the drone can see, and $\| \cdot \|$ represents the distance between the center of zone $z$ and the physical location of the drone.   When a set of observations are taken, the intensities of the fire in the set of observable zones are read using onboard sensors.  The set of observations at time $t+1$ is given by,
\bn
\label{eqn:obs_set}
\hat{\Ycal}_{t+1} = \{ \yhat_{t+1,z'}: z' \in \Zcal_{t+1}^{obs}\}.
\en
Additionally, we assume the intensity observations are very low noise and proportional to $Q_{t+1,z}$ through some intensity constant $c_z$; hence, $\yhat_{t+1,z} = c_z Q_{t+1,z} + \varepsilon_{t+1}$, where $\varepsilon_{t+1} \sim N(0, \sigma^2)$ with known sensor variance.

The random wind speed at time $t+1$, $\Uhat_{t+1}$, and the random wind direction, $\phihat_{t+1}$, are also imperfectly known to the controlling agent.  In this problem, we assume the wind speed and directions are drawn from stationary distributions over the time horizon.  The wind speeds at each time are drawn from an i.i.d. log-normal distribution with known parameters $\mu^U$ and $\sigma^{U}$, and the wind directions are drawn from a zero-mean Gaussian random-walk with known parameter, $\sigma^{\phi}$.  We assume the drone has on-board sensors which can measure the wind in real-time; hence, $U_t$ and $\phi_t$ are known perfectly at time $t$.

\paragraph{Belief State}
The belief state for the this belief model will contain the parameters of the dynamically changing probability distributions over each of the uncertain environment states.  In the fuel level model, there are three physical state variables for each zone which we do not know: the amount of healthy fuel ($H_{tz}$), the amount of burning fuel ($Q_{tz}$), and the amount of dead or extinguished fuel ($D_{tz}$).  The total level of fuel in each zone is constant; hence, the property in equation (\ref{eqn:fuel_prop}) can be used to construct a probability distribution.
\begin{lemma}
\label{def:1}
Assume the true fuel states $Q_{tz}$, $H_{tz}$, and $D_{tz}$ are discretized into integer valued groups constrained by $Q_{tz} + H_{tz} + D_{tz} = \eta_z$. Let $p^{Q}_{tz}$, $p^{H}_{tz}$, and $p^{D}_{tz}$ be the percentage of fuel in each of the respective states in the fuel level model.  Since, there is a fixed amount of fuel in each zone, these percentages will sum to 1.  Therefore, this structure induces a multinomial distribution on the state space given by,
\bn
Q_{tz}, H_{tz}, D_{tz} \sim Mult(\eta_z, p^{Q}_{tz}, p^{H}_{tz}, p^{D}_{tz}).
\en
\proof{proof.}
Follows directly from the definition of the multinomial distribution.
\endproof
\end{lemma}

The parameters of the multinomial distribution in lemma (\ref{def:1}) are unknown to the controller at time $t$, but can be estimated through partial observations.  The estimates for each of the parameters of the multinomial distributions will be denoted by $\pbar^{Q}_{tz}$, $\pbar^{H}_{tz}$, and $\pbar_{tz}^{D}$, respectively.  $K_{tz}$ is an indicator variable for whether zone $z$ has ignited at time $t$.  The belief state for $K_{tz}$ is a Bernoulli random variable with parameter $p_{tz}^K$.  Hence, the belief state is given by
\bn
\label{eqn:belief_state}
B_t = (\pbar_{tz}^K, \pbar_{tz}^Q, \pbar_{tz}^H, \pbar_{tz}^D)_{z\in\Zcal}.
\en
The initial belief state includes the prior beliefs about the dynamic parameter before any estimates, and the parameters of the stationary wind distributions.

\begin{figure}[ht]
\includegraphics[width=1\textwidth, center]{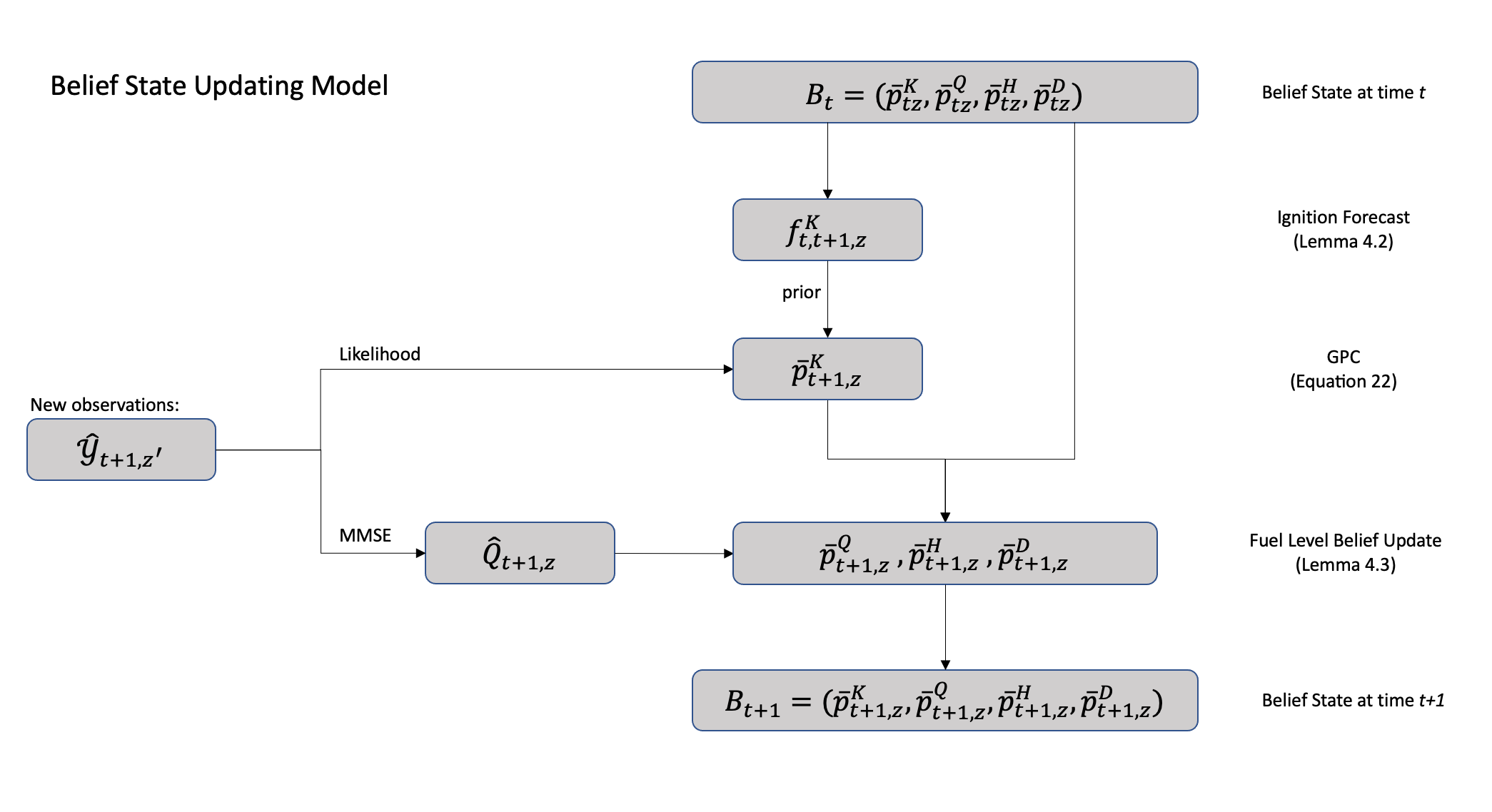}
\caption{Diagram depicting the steps in the updating procedure for the belief model. }
\label{fig:flowchart}
\end{figure}

\paragraph{Belief State Updates}
At time $t+1$, the drone will move to a new location, make a new set of observations, receive the helicopter decision and update the belief state with the new information it has.  Once the new observations are made, each belief state variable will need an updating equation.  First the parameters of the Bernoulli distributions will be updated, then the parameters of the multinomial distributions.  Figure \ref{fig:flowchart} shows a flowchart describing each step of the updating procedure.

The first step is to update the Bernoulli distributions in each zone.  Then, the updating procedure uses the assumed transition function to produce priors and the observations are used to construct a likelihood model.  The posterior is formed through a Bayesian update between the two estimates.   Lemma \ref{lem:forecast} shows how to forecast the prior distribution at $t+1$ and Equation \eqref{eqn:GPC_pred} shows how to predict the posterior by performing Gaussian process classification with the forecast model as the GP prior.

\begin{lemma}[Ignition Forecast]
\label{lem:forecast}
Assume we have a belief about whether each zone has ignited, $\pbar_{tz}^K$.  Before any observations, the prior knowledge we have about a zone $z$ at $t+1$ can be produced by a probabilistic forecast $f_{t,t+1,z}^K$. Let $\Lhat_{t+1,zz'}$ be the random length a fire will travel between time $t$ and $t+1$ in the direction of $z$ from zone $z'$.  Hence, the forecasting equations for $\pbar_{tz}^K$ are given by,
\bn
\label{eqn:postdecisionK}
\pbar_{tz}^{K,x} &=&  \pbar_{tz}^{K} \mathbf{1}_{[\Zcal \setminus \Zcal_{t}^{ext}]}  \\
\label{eqn:prior_drift}
f^K_{t,t+1,z} &=& 1 - \prod_{z' \in \Zcal} \left( 1 - \pbar_{tz'}^{K,x} + \Pbb[\Lhat_{t+1,zz'} < \|z - z'\| \big| K^x_{tz'} = 1] \pbar_{tz'}^{K,x}\right),
\en
where, equation \eqref{eqn:prior_drift} finds the probabilistic forecast that any surrounding zone will ignite zone $z$.
\proof{proof.}
See Appendix.
\endproof
\end{lemma}

Lemma \ref{lem:forecast} gives a sequence of equations for forecasting the belief of $K_{tz}$ for any zone.  However, the distribution of $\Lhat_{t+1,zz'}$ is impossible to know perfectly.  A common assumption in the wildfire literature is fire spreads elliptically as a function of the wind direction.  Hence, we propose the following corollary to Lemma \ref{lem:forecast}.

\begin{corollary}
Let $\phi_t$ be the wind direction at time $t$ measured by the drone.  Let $\psi_{zz'}$ be the angle between two zones $z$ and $z'$ with respect to north.  Let $\theta = (\theta_0, \theta_1) \in \mathbb{R}^2_+$ be parameters of an isotropic kernel function.  We assume the following CDF of $\Lhat_{t+1,zz'}$ is,
\bn
\label{eqn:prob_assume}
\Pbb[\Lhat_{t+1,zz'} < \|z - z'\| \big| K^x_{tz'} = 1] = 1 - \exp(-\gamma(\psi_{zz'}|\phi_t, \theta) \|z - z' \|),
\en
where $\gamma(\psi_{zz'}|\phi_t, \theta) = \theta_0 \cos(\psi - \phi_t) + \theta_1$ is a parameterized function of the interzone angles and wind direction.  Hence, the forecasts simplify to,
\bn
\label{eqn: prior_drift_param}
f^K_{t,t+1,z} = 1 - \prod_{z' \in \Zcal} \left( 1 - \pbar_{tz'}^{K,x}e^{-\lambda(\psi|\phi_t,\theta)\|z-z'\|} \right).
\en

\proof{proof.}
Substitute equation (\ref{eqn:prob_assume}) into equation (\ref{eqn:prior_drift}) and simplify.
\endproof
\end{corollary}
The forecasts will help support the statistical estimates of zones which have not been observed for a while.  The procedure for extracting the likelihood probabilities for each of the Bernoulli beliefs for $K_{tz}$ is outlined in the subsequent paragraphs.

Assume the new observations $\yhat_{t+1,z'} = Y_{t+1,z'} + \varepsilon_{t+1}$, where $\varepsilon_{t+1} \sim N(0, \sigma^2)$.  The measurement device noise, $\sigma^2$, is known.  Through equation \eqref{eqn:Kdef} we construct a hypothesis test to produce $\Khat_{t+1,z}$.  This produces a log-likelihood test given by,
\bn
\label{eqn:lltest}
\Khat_{t+1,z'} = \begin{cases}
      1 & \frac{\Phi\left(u\right)}{\psi\left(u\right)}> l , \\
      0 & otherwise.
   \end{cases}
\en
where $u = \frac{\yhat_{t+1,z'}}{\sigma}$ and $l$ is the risk adjusted classification threshold.  In this equation $\Phi(\cdot)$ is the normal CDF and $\phi(\cdot)$ is the normal PDF.

After transforming the intensity readings into boolean observations, the next step is to predict $K_{tz}$ in the unobserved zones, $z \in \Zcal \setminus \Zcal_{t+1}^{obs}$.  Assume our observations of $\Khat_{tz}$ are samples drawn from a function $\sigma(g(z))$, where $g(z)$ is a latent function with a Gaussian Process (GP) prior and $\sigma(\cdot)$ is a logistic function.  We can fully characterize the GP prior with a mean function  and a covariance kernel.  We define the prior mean over the latent function by passing the probability forecasts through the inverse logistic function,
\bns
m(z) = \sigma^{-1}(f^K_{t,t+1,z}),
\ens
and an RBF covariance kernel parameterized by $\theta^{cov}$, given by,
\bn
\label{eqn:covariance}
Cov[g(z), g(z')|\theta^{cov}] = \theta^{cov}_0 exp(-\theta_1^{cov} \|z - z'\|^2) + \theta_2^{cov} \mathbf{1}_{[z = z']}  \hspace{4mm}  \forall z,z' \in \Zcal.
\en
For any unobserved zone $z$, the approximate posterior mean and variance of the latent function can be computed with the Laplace approximation (\cite{GaussianProcess}).  Let the predicted mean of the posterior latent function be given by $\bar{g}(z)$ and the corresponding latent predictive variance be $\mathbb{V}[g(z)|\hat{\Ycal}_{t+1}]$.  The probabilistic prediction for the unobserved zones igniting is computed via the method in \cite{MacKay}.  Hence, $\pbar_{t+1,z}^K$ is given by,
\bn
\label{eqn:GPC_pred}
\pbar_{t+1,z}^K = \begin{cases}
      \sigma\left(\kappa(g(z)|\hat{\Ycal}_{t+1}) \bar{g}(z)\right) & z \in \Zcal \setminus \Zcal_{t+1}^{obs}, \\
      \Khat_{t+1,z} & z \in \Zcal_{t+1}^{obs}.

   \end{cases}
\en
where $\kappa\left(g(z)|\hat{\Ycal}_{t+1}\right) = \sqrt{(1 + \frac{\pi}{8} \mathbb{V}[g(z)|\hat{\Ycal}_{t+1}])^{-1}}$.



The updated belief $\pbar_{t+1,z}^K$ is used to update the belief about the fuel levels in each zone.  Lemma \ref{lem:fuel_updates} provides the updating equations for the belief about the fuel levels in each zone.

\begin{lemma} [Fuel Level Belief Update]
\label{lem:fuel_updates}
The belief state updates for the fuel levels will be different for observed versus unobserved zones.  For an observed zone $z\in \Zcal_{t+1}^{obs}$, let $c_z$ be the known ignition constant between zones.  The intensity observations can be used to compute the minimum mean square estimate (MMSE) of $Q_{tz}$ using the observations through $\Qhat_{t+1,z} = \frac{\yhat_{t+1,z}}{c_z}$.  The observed zones have low noise, hence the MMSE estimate of $Q_{t+1,z}$ is used as the new belief state.  To mitigate risk, we assume the remaining fuel in the zone is susceptible to being burned and there is no dead fuel yet.

For unobserved zones, the observations are used in conjunction with the transition equations we assume about the environment agent.  Let $\pbar_{t+1,z}^{start}$ be the probability that fire has been started for the first time in zone $z$ at time $t+1$.   Let $\pbar_{t+1,z}^{burn}$ be the probability that fire is continuing to burn from the previous time step in zone $z$ at time $t+1$.  We approximate the terms by the following expressions,
\bns
\pbar_{t+1,z}^{start} &=& \pbar_{t+1,z}^K -  \frac{\pbar_{t+1,z}^K \pbar_{tz}^{K,x}}{f_{t,t+1,z}},\\
\pbar_{t+1,z}^{start} &=& \frac{\pbar_{t+1,z}^K \pbar_{tz}^{K,x}}{f_{t,t+1,z}}.
\ens
Note, if $f_{t,t+1,z} = 0$, then we define $\pbar_{t+1,z}^{start} = \pbar_{t+1,z}^K$ and $\pbar_{t+1,z}^{burn} = 0$.

The belief about the fuel is updated by taking the conditional expectation of equations \eqref{eqn:env_Qt1_update} - \eqref{eqn:env_Dt1_update} with respect to the multinomial belief state at $t$ and the set of Bernoulli belief states at $t+1$.  Hence,
\bn
\label{eqn:Qt1_Eupdate}
\pbar_{t+1,z}^Q &=& \E_t\left[\frac{Q_{t+1,z}}{\eta_z}\right] = \begin{cases}
      \left(1 - \xi_z + \lambda_z (\eta_z - 1) \pbar_{tz}^H \right) \pbar_{tz}^Q \pbar_{t+1,z}^{burn} + \frac{Q^{init}_{z}}{\eta_z} \pbar_{t+1,z}^{start} & z \in \Zcal \setminus \Zcal_{t+1}^{obs}, \\
      \frac{\Qhat_{t+1,z}}{\eta_z} & z \in \Zcal_{t+1}^{obs},
\end{cases} , \\
\label{eqn:Ht1_Eupdate}
\pbar_{t+1,z}^H &=& \E_t\left[\frac{H_{t+1,z}}{\eta_z}\right] =  \begin{cases}
      \pbar_{tz}^H - \lambda_z (\eta_z-1)\pbar_{tz}^Q\pbar_{tz}^H \pbar_{t+1,z}^{burn} - Q^{init}_{z} \pbar_{t+1,z}^{start} & \hspace{8.5mm} z \in \Zcal \setminus \Zcal_{t+1}^{obs}, \\
       \frac{\eta_z - \Qhat_{t+1,z}}{\eta_z} & \hspace{8.5mm} z \in \Zcal_{t+1}^{obs},
\end{cases}, \\
\label{eqn:Dt1_Eupdate}
\pbar_{t+1,z}^D &=& \E_t\left[\frac{D_{t+1,z}}{\eta_z}\right] = \begin{cases}
      \pbar_{tz}^D + \pbar_{tz}^{Q} (1 - \pbar_{t+1,z}^{burn}) + \xi_z \pbar_{tz}^{Q} \pbar_{t+1,z}^{burn} & \hspace{22mm} z \in \Zcal \setminus \Zcal_{t+1}^{obs}, \\
      0 & \hspace{22mm} z \in \Zcal_{t+1}^{obs},
\end{cases}
\en
where $\eta_z$, $\lambda_z$, and $Q^{init}_z$ are known fuel constants for each zone.
\proof{proof.}
See Appendix.
\endproof
\end{lemma}

To summarize the belief state update for each zone $z \in \Zcal$:
\begin{enumerate}
    \item Forecast the belief about $K_{t+1,z}$ (Equation \eqref{eqn: prior_drift_param}).
    \item Use GPC to extract $\pbar_{t+1,z}^K$ from observations and blend them with the GP prior (Equation \eqref{eqn:GPC_pred}).
    \item Update the fuel level belief states (Equations \eqref{eqn:Qt1_Eupdate} - \eqref{eqn:Dt1_Eupdate}).
\end{enumerate}

\subsubsection{Drone Model}
\label{sec:drone_model}
The drone model will contain the five components of the unified framework: the drone state variable, drone decision variable, the drone exogenous information process, the drone transition function, and the objective function.

\paragraph{Drone State Variable}  The drone state variable consists of all information needed to update the transition function, compute the policy, and evaluate the objective function.  The pre-decision state variable is given by,
\bn
S_t^{drone} = \left(R_t^{drone}, R_t^{batt}, U_t, \phi_t, \Zcal_{t-1}^{ext}, B_t\right).
\en
where
\bns
R_t^{drone} &=& (R_{tz}^{drone})_{z\in\Zcal} = \text{binary vector representing the physical location of the drone},\\
R_t^{batt} &=& \text{battery level of drone},\\
U_t &=& \text{wind speed at time $t$},\\
\phi_t &=& \text{wind direction at time $t$},\\
\Zcal_{t-1}^{ext} &=& \text{set of extinguished zones through time $t-1$},\\
B_t &=& \text{belief state at time $t$}.
\ens

The initial drone state variable holds the initial states of each dynamic state variable as well as any static information as latent state variables, given by
\bn
S_0^{drone} = \left(R_0^{drone}, R_0^{batt}, U_0, \phi_0, B_0, \Xi_z \right).
\en

\paragraph{Drone Decision Variable}
The drone must decide which zone to travel to at time $t$ to make new observations at $t+1$.  The decision variable for each zone is given by,
\bns
x_{tz}^{drone} &=& \begin{cases}
      1 & \text{drone decides to move to zone $z$}, \\
      0 & otherwise,
   \end{cases}
\ens
cosntrained by $\sum_{z\in\Zcal} x_{tz}^{drone} = 1$.   The drone has a maximum speed which it can travel, hence it has a maximum distance it can travel in one time step, $d^{max}$.  The binary vector representing the entire decision variable is given by, $\xbf_{t}^{drone} = (x_{tz}^{drone})_{z\in\Zcal}$.

\paragraph{Drone Exogenous Information} The drone exogenous information is the set of information which arrives between time $t$ and $t+1$.  The observations at $t+1$, the wind speed at $t+1$, the wind direction at $t+1$, and the decision from the other agents.  The exogenous information is given by, $W_{t+1}^{drone} = (\hat{\Ycal}_{t+1}, \Uhat_{t+1}, \phihat_{t+1}, \xbf_{t}^{heli})$.

\paragraph{Drone Transition Function} The transition function is the set of equations describing how each state variable is updated.  The physical state at $t+1$ is the decision from $t$; hence, $R_{t+1,z}^{drone} = x_{tz}^{drone}$.  The wind speed and direction evolve exogenously.  The set of extinguished zones is updated through equation \eqref{eqn:env_ext_update}.  The multinomial belief state parameters are updated through equations \eqref{eqn:Qt1_Eupdate} - \eqref{eqn:Dt1_Eupdate}.  The Bernoulli belief state over $K_{t+1,z}$ is updated with equation \eqref{eqn:GPC_pred}.   The battery level depletes at a constant rate; hence,
\bn
\label{eqn:batt_trans}
R_{t+1}^{batt} = \max\{R_t^{batt} - 10, 0\}.
\en

\paragraph{Drone Objective Function} The joint goal of the collaborative agents is to minimize the cumulative cost of fuel burned.  At each time step, the one-step contribution function is given by,
\bn
\label{eqn:one_step_cost}
C(S_t^{env}, \xbf_t^{heli}) = \sum_{z\in\Zcal} r_z Q_{t+1,z} + c^{fail} \mathbb{I}_t,
\en
where $c^{fail}$ is a penalty cost for the drone not returning home and
\bn
\label{eqn:indicator_batt}
\mathbb{I}_t = \begin{cases}
    1 & \text{if $R_t^{batt} = 0$ and $R_{tz}^{drone}$ = 1 for any $z \in \Zcal \setminus z_0$ } \\
    0 & \text{otherwise}.
\end{cases}
\en
In other words, the agent is penalized for not returning home before it runs out of battery.  The helicopter decision will directly impact the environment which is why it is the argument in the one-step contribution function.  However, the policy for the helicopter is a function of the belief state which is updated through observations chosen by the drone policy.  Hence, the drone objective is given by the following optimization problem,
\bn
\label{eqn:drone_objective}
\min_{\pi^{drone} \in \Pi^{drone}} \E \left\{ \sum_{t=0}^{T-1} C(S_t^{env}, X^{\pi^{heli}}(S_t^{heli}(X^{\pi^{drone}}))) | S_0^{drone} \right\},
\en
where the second argument is the helicopter policy which is a function of the helicopter state.  The helicopter state is determined by the belief state communicated by the drone.  Hence, the helicopter decision is a function of the drone policy.

\subsection{Helicopter Agent}
\label{sec:heli_model}
The helicopter model will contain the five components of the unified framework.

\paragraph{Helicopter State Variable}  The helicopter state variable consists of all information needed to update the transition function, compute the policy, and evaluate the objective function.  The pre-decision state variable is given by, $S_t^{heli} = (B_t)$.

\paragraph{Helicopter Decision Variable} The decision variable describes the decision which the helicopter makes at each time step along with any constraints.  The decision variable is defined by,
\bns
x_{tz}^{heli} = \begin{cases}
      1 & \text{extinguish zone $z$}, \\
      0 & otherwise,
   \end{cases}
\ens
constrained by $\sum_{z\in\Zcal} x_{tz}^{heli} = 1$.  The binary vector of helicopter decisions is given by, $\xbf_t^{heli}$.

\paragraph{Helicopter Exogenous Information} The exogenous information describes the information that is first available to the helicopter between time $t$ and $t+1$.  The exogenous information variable is the new belief state, $W_{t+1}^{heli} = B_{t+1}$.

\paragraph{Helicopter Transition Function} The transition function describes the equations to update the state variable.  The state variable of the helicopter is the information about the belief state which is communicated by the drone. Thus, the helicopter state at time $t+1$ is $S_{t+1}^{heli} = B_{t+1}^{drone}$.

\paragraph{Helicopter Objective function}
The helicopter objective function describes the performance metric used to create the policy.  The joint goal of the collaborative controllers is to mitigate the damage caused by burning fuel. The one-step cost function for the helicopter is given by equation \eqref{eqn:one_step_cost}.  The search over policies for the helicopter is given by the following optimization problem,
\bn
\label{eqn:heli_objective}
\min_{\pi^{heli} \in \Pi^{heli}} \E \left\{ \sum_{t=0}^{T-1} C(S_t^{env}, X^{\pi^{heli}}(S_t^{heli}) | S_0^{heli} \right\}.
\en.

\section{Designing Policies}
\label{sec:policies}
To solve a sequential decision problem, we must design a policy which outputs decisions based on the state of the system.  The policy, $\pi$ is a function which takes the state of the agent as input and maps it into a decision $X^{\pi}$.  In theory, the optimal policy is a mapping which makes the best decision for any state with respect to the sum of expected rewards for the remaining time horizon.  However, due to the curse of dimensionality (\cite{PowellEJOR}) it is not possible to find the feasible policy when the state space becomes too large for a fully observable problem.  In this wildfire problem formulation, we have partial observability with large state spaces, decision spaces, and observation spaces; hence, there is no feasible way to find the optimal policy.  When there is no feasible way to achieve optimality, we must search over the set of approximate policies for the best solution.

There are two strategies for designing policies: policy search and lookahead approximations.  The two strategies can be further subdivided into four classes of policies encompassing the set of all possible strategies.  The policy search class contains the policy function approximation (PFA) class of policies and the cost function approximation (CFA) class of policies.  The lookahead approximation strategy contains the value function approximation (VFA) class of policies and the direct lookahead approximation (DLA) class of policy.  The four classes of policies are defined as:

•	\textbf{Policy Function Approximations (PFA)} are analytical functions which directly map states into decisions.  These functions usually contain hyperparameters which must be tuned by optimizing a simulator or online learning in the real world.  Examples of PFAs are linear parametric functions, nonlinear parameteric functions, or look-up tables,

•	\textbf{Cost Function Approximations (CFA)} are optimization models which contain hyperparameters (usually designed to account for uncertainty).  The CFA requires hyperparameter tuning within a simulator or online in the real world.  Examples of CFAs are buffer stock when optimizing supply chains, slack scheduling when optimizing airline fleets, or upper confidence bounding in bandit problems,

•	\textbf{Value Function Approximations (VFA)} are lookahead approximations which seek to approximate the value function in the optimal policy and then solve the approximate version of Bellman’s optimality equation.  Examples of VFAs can be found in the fields of approximate dynamic programming (\cite{Powell_ADP}), reinforcement learning (\cite{sutton2018reinforcement}), and SDDP (\cite{birge2011introduction}),

•	\textbf{Direct Lookahead Approximations (DLA)} are lookahead strategies which optimize models which look directly into the future using some approximations.  Examples of DLAs are model predictive control (\cite{agachi20162}), and Monte Carlo Tree Search (\cite{browne2012survey}).

It is also possible to solve hybrid policies between the four classes.  An example of a hybrid model would be parameterizing a DLA with tunable parameters to adjust for risks in the future.

We must design a policy for each agent to make the best decisions with respect to its goal.  The helicopter policy will be the best solution to equation \eqref{eqn:heli_objective}, and the drone policy will be the best solution to equation \eqref{eqn:drone_objective}.  Subsection \ref{sec:heli_policies} will present two policies for optimizing helicopter decisions and subsection  \ref{sec:drone_policies} will present three policies for optimizing drone decisions.

\subsection{Helicopter Policies}
\label{sec:heli_policies}
In this section, we design two policies to extinguish the fire with the helicopter agent.  The first policy predicts where the most expensive fires will start in the next time step and chooses to extinguish the zone with the highest probability of extinguishing.  The second policy is a CFA-DLA hybrid policy which forecasts where fire will start at $t+1$ and trades off this forecast with a tunable parameter for extinguishing the highest cost of burning fuel now.  This will preemptively extinguishes the fire-front where we expect it to have a high cost of burning fuel and will spread in the future.

\subsubsection{Helicopter probabilistic DLA}
The probabilistic DLA policy (DLA-1) for the helicopter agent will use the belief state at time $t$ to estimate the probability of burning fuel in each zone in the future.  Then a filter with the shape of the extinguishing region is passed over the entire region and the zone with the probability of extinguishing the most expensive new fire is chosen.  The probabilistic DLA policy is given by,
\bns
X^{DLA-1}(S_t^{heli}) = \max_{\xbf_{t}^{heli} \in \Xcal_{t}^{heli}} \sum_{z\in \Zcal_{t}^{ext}(\xbf_t^{heli})} r_z (f_{t,t+1z}^K - \pbar_{tz}^K).
\ens

\subsubsection{Helicopter CFA-DLA hybrid}
The goal of this policy is to forecast which zones will ignite in the future and trade this computation off with the cost of burning fuel at time $t$.  We can preemptively decide where to extinguish now based on where we believe the most costly burning fuel is, and then weight this by the probability it will ignite in the future.  Then, a filter is passed over the region to produce a map of the CFA-DLA cost heuristic for each decision.  The helicopter policy then will choose the location with the largest decrease in expected future fires extinguished.  The policy is given by,
\bns
X^{CFA-DLA}(S_t^{heli}) &=& \argmax_{\xbf_t^{heli}} \E_t\left[ \sum_{z\in\Zcal_t^{heli}(\xbf_t^{heli})} r_z Q_{tz}^x + \theta^{heli} (K_{t+1,z} - K_{tz}) \right], \\
&=& \argmax_{\xbf_t^{heli}}\sum_{z\in\Zcal_t^{heli}(\xbf_t^{heli})} r_z \eta_z \pbar_{t+1,z}^{Q} + \theta^{heli}f_{t,t+1,z}^{start}, \\
\ens
where $\theta^{heli}$ is a tunable parameter used to trade-off between the two terms.  Through inspection of equation \eqref{eqn:env_Qtx_update} we can see this is a heuristic way of predicting the next burning fuel level without knowing the true dynamics of the environment.

\subsection{Drone Policies}
\label{sec:drone_policies}
The drone policies are designed to collect the most valuable information to guide the helicopter, but they are constrained by the requirement to return home at the end of the time horizon.  The drone policies will require logic which balances the goal of collecting information, but also evaluating whether it needs to begin planning a path home.

We present multiple policies motivated from the Bayesian optimization literature to solve the problem.  Bayesian optimization (BO) usually seeks solutions with less than 25 dimensions, however we are able to adapt some of the policies from the community to overcome some of the limitations of our large action space.  However, the size of the action space will limit the set of policies we can consider.

\subsubsection{Myopic PFA-CFA hybrid}
We design a drone policy which is a hybrid policy between a PFA and a CFA.  Recall, a PFA is an analytical function which maps states directly into actions, and usually contains policies which use simple rules.  In this policy, we use a simple rule to decide when the drone needs to return home by checking whether there is enough battery life to get home in the next time step.  The PFA portion of the policy creates the decision whether to continue looking for new observations or turn back and go home.  Let $\xbf_{tz_0}^{drone}$ be a shorthand notation for the decision to return home.  The PFA output is a Boolean variable whether to return home or choose a zone to observe with the CFA.  The PFA is given by,
\bns
X^{PFA}(S_t^{drone}) = \begin{cases}
    1 & \text{if $\|z_0 - R_t^{drone}\| \leq d^{max} * (R_{t}^{batt} - \Delta^{batt} * (t+1))$},\\
    0 & \text{otherwise}.
\end{cases}
\ens

The CFA portion of the policy decides which zones are the best to travel to using an interval estimation metric to capture the value of making a new observation.

The interval estimation metric is used to determine the marginal value of observing the intensity in a single zone.  The interval estimation metric for zone $z$ is given by,
\bns
m_{tz} = \mu_{tz} + \theta^{IE} \sigma_{tz},
\ens
where $\theta^{IE}$ is a tunable parameter, and
\bns
\mu_{tz} &=& \E[r_z Q_{tz}] = r_z \eta_z \pbar_{tz}^Q, \\
\sigma_{tz} &=& r_z \sqrt{\E[Q_{tz}^2] - \E[Q_{tz}]^2} = r_z \sqrt{\eta_z \pbar_{tz}^Q (1-\pbar_{tz}^Q)}.
\ens

The CFA policy performs a sum over the image masked with the zones within the observable set for each possible decision.  Hence, the CFA portion of the policy becomes,
\bns
X^{CFA}(S_{t}^{drone}) = \argmax_{\xbf_{t}^{drone} \in \{0,1\}^{|\Zcal|}} \sum_{z'\in\Zcal^{obs}_{t+1}(\xbf_{t}^{drone})} m_{tz'}.
\ens

If the decision chosen by the CFA policy is not in the feasible decision set, then we project it back to the closest point in the feasible decision set.  Let $\Pi_{\Xcal}$ be the projection operator which projects $\xbf_t^{drone}$ back onto closest point in the feasible decision set, $\Xcal_{tz}^{drone}$.  Hence, the hybrid policy becomes,
\bn
X^{hybrid}(S_t^{drone}) = \Pi_{\Xcal}\left[X^{CFA}(S_t^{drone}) \left(1-X^{PFA}(S_t^{drone})\right) + \xbf_{tz_0}^{drone} X^{PFA}(S_t^{drone})\right].
\en

\subsubsection{Multi-stage Thompson sampling lookahead approximation}
In a multi-stage Thompson sampling lookahead approximation (TS-DLA), we must generate an approximate model to be used to simulate the future.  We will use the five components of the unified framework to simplify the drone model from section \ref{sec:drone_model}.  The general lookahead variable notation is accented with a tilde and includes two time indices.  The first time index represents the time it is in the base model and the second time index represents the time it is in the lookahead model.

\begin{itemize}
    \item \textbf{Lookahead state variable:} We approximate the state variable in this lookahead model by sampling the belief state at time $t$, then generating multiple scenarios into the future to optimize.  The general state variable starting at time $t$ in the base model and looking forward to a future time $t'$ is given by,
    \bns
    \Stilde_{tt'}^{drone} = \left(\tilde{R}_{tt'z}^{drone}, \tilde{R}_{tt'z}^{batt}, \tilde{K}_{tt'z}, \tilde{Q}_{tt'z}, \tilde{H}_{tt'z}, \tilde{D}_{tt'z}, \tilde{U}_{tt'}, \tilde{\phi}_{tt'}, \tilde{\Zcal}\right)_{z\in\Zcal},
    \ens
    where, $\tilde{R}_{tt'z}^{drone}$, $\tilde{R}_{tt'z}^{batt}$, $\tilde{\Zcal}_{tt'}$, $\tilde{U}_{tt'}$, $\tilde{\phi}_{tt'}$, are known deterministically at time $t$ within the lookahead model.   The remaining lookahead state variables are sampled through the procedure in Appendix \ref{app:sampling} to generate $\Stilde_{tt}^{drone}$.

    \item \textbf{Lookahead decision variable:} The lookahead decision variable is equivalent to the drone decision variable from section \ref{sec:drone_model}.

    \item \textbf{Lookahead exogenous information:} The lookahead exogenous information includes the new wind speed and direction, the binary adjacency matrix from the environment transition function, and the helicopter decisions.  The exogenous information variable is given by,
    \bns
    W_{t,t+1}^{drone} = \left(\tilde{U}_{t,t'+1}, \tilde{\phi}_{t,t'+1}, \tilde{\mathbf{A}}_{t,t'+1}, \tilde{\xbf}_{t,t'+1}^{heli}\right),
    \ens
    where $\tilde{U}_{t,t'+1}$ and $\tilde{\phi}_{t,t'+1}$ are sampled from the distributions assumed in section \ref{sec:belief_model}, $\tilde{\mathbf{A}}_{t,t'+1}$ is sampled through the approximate surface spread model described in Appendix \ref{app:Rothermel}, and $\tilde{\xbf}_{t,t'+1}^{heli}$ is computed by passing the lookahead state variable through the helicopter policy which is communicated to the drone.

    \item \textbf{Lookahead transition function:} The lookahead transition function for each of the sampled belief state variables are given by equations \eqref{eqn:env_ext_update} - \eqref{eqn:env_Dt1_update} with the samples as inputs.  The lookahead transition for the new drone position is the decision from the previous step.  The battery levels transition according to equation \eqref{eqn:batt_trans}.

    \item \textbf{Lookahead objective function:} The lookahead objective function is different than the objective function from the base model in section \ref{sec:drone_model} because we form an approximate one step contribution function to promote learning by maximizing the amount of information we collect.  The approximate cost function is given by,
    \bns
    \tilde{C}\left(\Stilde_{tt'}^{drone}, \tilde{\xbf}_{tt'}^{drone}\right) = \sum_{z\in\Zcal_{tt'}^{obs}( \tilde{\xbf}_{tt'z}^{drone})} r_z \left(\tilde{Q}_{t,t'+1,z}\right) - c^{fail} \mathbb{I}_{tt'}.
    \ens
    where $\mathbb{I}_t$ is defined in equation \eqref{eqn:indicator_batt}.
\end{itemize}

We solve the approximate model by choosing a truncated horizon length $H$ and simulating $m$ scenarios of the lookahead model $H$ steps into the future.  The objective function will decompose into the following optimization problem,
\bn
\label{eqn:sampled_opt_problem}
\max_{\tilde{\xbf}_{tt}^{drone}\in \tilde{\Xcal}_{tt}^{drone}, ... , \tilde{\xbf}_{t,t+H}^{drone}\in\tilde{\Xcal}_{t,t+H}^{drone}} \tilde{C}\left(\Stilde_{tt}^{drone}, \tilde{\xbf}_{tt}^{drone}\right) + \frac{1}{m} \sum_{i=1}^{m} \sum_{t'=t+1}^{t+H} \tilde{C}\left(\Stilde_{tt'}^{drone}(\tilde{\omega}^{i}), \tilde{\xbf}_{tt'}^{drone}\right),
\en
where $(\tilde{\omega}^i)_{i=1}^{m} \subset (\Wcal_{tt'}^{drone})_{t'=t+1}^{t+1+H}$.  In other words, we generate a set of sample paths through the set of possible exogenous information sets.  Hence, $\Stilde_{tt'}^{drone}(\tilde{\omega}^{i})$ is the future state at $t'$ generated by simulating the lookahead variable with information path, $\tilde{\omega}^i$.

\begin{lemma}
\label{lem:samp_LA}
Let $\tilde{\xbf}_{tt}^{drone}$, $\tilde{\xbf}_{t,t+1}^{drone}$, ... , $\tilde{\xbf}_{t,t+H}^{drone}$ be the binary vectors representing the drone decisions $H$ steps into the lookahead model.  Let $\Mbf \in \{0,1\}^{|\Zcal| \times |\Zcal|}$ be a matrix where each column is a binary mask over the field of view for each possible location in the set of zones.  Let $\Vbf \in \{0,1\}^{|\Zcal| \times |\Zcal|}$ be a matrix where each column is a binary mask over the zones reachable within 10 minutes for each possible location in the set of zones.  Let $\tilde{\cbf}_{tt'} = \left[\frac{1}{m} \sum_{i=1}^{m} r_z (\tilde{Q}_{t,t'+1,z}(\tilde{\omega}^{i}) - \tilde{Q}_{tt'z}^x(\tilde{\omega}^{i}))\right]_{z\in\Zcal} - c^{fail} \mathbb{I}_{tt'}$ be a vector where each entry represents the approximate marginal cost for each individual zone if it is observed at $t'$.  Then, the policy induced by the optimization problem in equation \eqref{eqn:sampled_opt_problem} reduces to an integer linear program given by,
\bn
\label{eqn:samp_LA}
\tilde{\xbf}_{tt}^{drone, *}, ... , \tilde{\xbf}_{t,t+H}^{drone, *} &=& \arg\max_{\tilde{\xbf}_{tt}^{drone}, ... , \tilde{\xbf}_{t,t+H}^{drone}} \sum_{t'=t}^{t+H} \tilde{\cbf}_{tt'}^T \Mbf \tilde{\xbf}_{tt'z}^{drone}, \\
\label{eqn:linint_constr1s}
\text{subject to} &:& \sum_{z\in\Zcal} \tilde{\xbf}_{tt'z}^{drone} = 1  \hspace{4mm} \forall t'\in{t,...,t+H}, \\
\label{eqn:linint_constr2s}
& & \Vbf \tilde{\xbf}_{tt'} -  \tilde{\xbf}_{t,t'+1} \geq 0 \hspace{4mm} \forall t'\in{t+1,...,t+H-1}, \\
\label{eqn:linint_constr3s}
& & \Vbf \tilde{R}_{tt} -  \tilde{\xbf}_{t,t+1} \geq 0 \hspace{4mm}, \\
\label{eqn:linint_constr4s}
& & \tilde{\xbf}_{tt}^{drone}, ... , \tilde{\xbf}_{t,t+H}^{drone} \in \{0,1\}^{|\Zcal|} \hspace{4mm} \forall t'\in{t,...,t+H-1},
\en
where constraint \eqref{eqn:linint_constr1s} ensures at each time step the decision chooses a single location, constraints \eqref{eqn:linint_constr2s} and \eqref{eqn:linint_constr3s} ensure that it can travel to the next location in a single time step, and constraint \eqref{eqn:linint_constr4s} keeps the variables binary valued.  This optimization problem has $H|\Zcal|$ variables and can be solved by any software package with a integer linear program solver.  The decision at time $t$ is given by $X^{TS-DLA}(S_t^{drone}) =\tilde{\xbf}_{tt}^{drone, *} $.
\proof{proof.}
See Appendix.
\endproof
\end{lemma}
This policy is capable of choosing observations by simulating the future environment and helicopter decisions in order to choose observations at time $t$ which will be the most helpful to the helicopter, but will likely not be extinguished in the future.  Furthermore, the approximate state variable at time $t$ helps balance exploration versus exploitation in the lookahead in a similar paradigm in which Thompson sampling does because we are sampling from the belief state then simulating the future (see \cite{shahriari2015taking}).

\subsubsection{Multi-stage Interval Estimation lookahead approximation}
In the multi-stage interval estimation lookahead (IE-DLA), we create our approximate model of the future by simulating the probability distributions over the uncertain state variables into the future opposed to sampling them.  Then, we use an interval estimation metric with a set of tunable parameters, $\theta^{IE},..., \theta_{t+H}^{IE}$, to promote learning in the future.  Each step in the future will have a different parameter because we will learn the tradeoff between exploitation versus exploration as a function of how far we are looking ahead.

We will use the five components of the unified framework to simplify the drone model from section \ref{sec:drone_model}.  The general lookahead variable notation is accented with a tilde and includes two time indices.  The first time index represents the time it is in the base model and the second time index represents the time it is in the lookahead model.  Furthermore, any expectation using a belief state distribution in the lookahead model is given by $\E[\cdot|\tilde{B}_{tt'}] = \Etilde_{tt'}[\cdot]$.

\begin{itemize}
    \item \textbf{Lookahead state variable:} We approximate the state variable in this lookahead model by sampling the belief state at time $t$, then generating multiple scenarios into the future to optimize.  The general state variable starting at time $t$ in the base model and looking forward to a future time $t'$ is given by,
    \bns
    \Stilde_{tt'}^{drone} = \left(\tilde{R}_{tt'z}^{drone}, \tilde{R}_{tt'}^{batt},  \tilde{\phi}_{tt'}, \tilde{\Zcal}_{tt'},  \Btilde_{tt'}\right)_{z\in\Zcal}.
    \ens
    \item \textbf{Lookahead decision variable:} The lookahead decision variable is equivalent to the drone decision variable from section \ref{sec:drone_model}.

    \item \textbf{Lookahead exogenous information:} The lookahead exogenous information includes the new wind speed and direction, the binary adjacency matrix from the environment transition function, and the helicopter decisions.  The exogenous information variable is given by,
    \bns
    W_{t,t+1}^{drone} = \left(\tilde{\phi}_{t,t'+1}, \tilde{\xbf}_{t,t'+1}^{heli}\right),
    \ens
    where $\tilde{\phi}_{t,t'+1}$ is sampled from the distributions assumed in section \ref{sec:belief_model}, and $\tilde{\xbf}_{t,t'+1}^{heli}$ is computed by passing the lookahead state variables through the helicopter policy which is known to the drone.

    \item \textbf{Lookahead transition function:} The lookahead transition functions are computed by forecasting each of the belief state variables.  The indicator variables which simulate the spread of the fire are computed with equation \eqref{eqn: prior_drift_param}.  The fuel state forecasts are computed with equations \eqref{eqn:Qt1_Eupdate} - \eqref{eqn:Dt1_Eupdate} and we use the equations for unobserved zones.  The lookahead transition for the new drone position is the decision from the previous step.  The battery levels transition according to equation \eqref{eqn:batt_trans}.

    \item \textbf{Lookahead objective function:} The lookahead objective function is equivalent to the objective function from the base model in section \ref{sec:drone_model}, but we form an approximate one step contribution function to promote learning.  The one-step contribution function is given by,
    \bns
    \tilde{C}\left(\Stilde_{tt'}^{drone}, \tilde{\xbf}_{tt'}^{drone}\right) = \sum_{z\in\Zcal_{tt'}^{obs}( \tilde{\xbf}_{tt'z}^{drone})} \tilde{\mu}_{tt'z} + \theta^{IE} \tilde{\sigma}_{tt'z},
    \ens
    where $\tilde{\mu}_{tt'z} = \Etilde_{tt'}\left[r_z \left(\tilde{Q}_{t,t'+1,z} - \tilde{Q}_{ttz}^x\right)\right]$, $\tilde{\sigma}_{tt'z} = r_z \sqrt{\left(\Etilde_{tt'}\left[ \left(\tilde{Q}_{t,t'+1,z} - \tilde{Q}_{ttz}^x\right)^2\right] - \tilde{\mu}_{tt'z}^2\right)}$, and $\theta^{IE}$ is a tunable parameter.
\end{itemize}

We solve the approximate model by choosing a truncated horizon length $H$ and using an marginal interval estimation metric from the Bayesian optimization literature to promote learning within in the lookahead.  We use the term marginal interval estimation as the one-step contribution.  The vector of marginal interval estimation contributions for each zone is given by $\tilde{\cbf}_{tt'}$.   We omit the $z$ subscript to denote the vector of means and variances given by $\tilde{\mu}_{tt'}$ and $\tilde{\sigma}_{tt'}$, respectively.  Furthermore, we use the approximation that an observed zone in the lookahead has zero variance for the remainder of the lookahead time horizon.

The objective function will decompose into the following optimization problem,
\bn
\label{eqn:prob_opt_problem}
\max_{\tilde{\xbf}_{tt}^{drone}\in \tilde{\Xcal}_{tt}^{drone}, ... , \tilde{\xbf}_{t,t+H}^{drone}\in\tilde{\Xcal}_{t,t+H}^{drone}} \tilde{\cbf}_{ttz}^T\Mbf \xbf_{tt}^{drone}  + \sum_{t'=t+1}^{t+H} \tilde{\cbf}_{tt'z}^T\Mbf \xbf_{tt'}^{drone} - (\xbf_{t,t'-1}^{drone})^T\Mbf^T diag(\tilde{\sigma}_{tt'})\Mbf\xbf_{tt'}^{drone},
\en
where $diag(\tilde{\sigma}_{tt'})$ represents the vector of standard deviations placed along the diagonal of a matrix of appropriate size.  Essentially, we set the standard deviation from the previously observed zones to zero.

\begin{lemma}
\label{lem:prob_LA}
Let $\tilde{\xbf}_{tt}^{drone}$, $\tilde{\xbf}_{t,t+1}^{drone}$, ... , $\tilde{\xbf}_{t,t+H}^{drone}$ be the binary vectors representing the drone decisions $H$ steps into the lookahead model.  Let $\Mbf \in \{0,1\}^{|\Zcal| \times |\Zcal|}$ be a matrix where each column is a binary mask over the field of view for each possible location in the set of zones.  Let $\Vbf \in \{0,1\}^{|\Zcal| \times |\Zcal|}$ be a matrix where each column is a binary mask over the zones reachable within 10 minutes for each possible location in the set of zones.  Let $\tilde{\cbf}_{tt'} = \left[\mu_{tt'z} + \theta^{IE} \sigma_{tt'} \right]_{z\in\Zcal}$ be a vector where each entry represents the approximate marginal interval estimation contribution for each individual zone if it is observed at $t'$.  Then, the policy induced by the optimization problem in equation \eqref{eqn:prob_opt_problem} reduces to an quadratic integer program with linear constraints given by,
\bn
\label{eqn:prob_LA}
\tilde{\xbf}_{tt}^{drone, *}, ... , \tilde{\xbf}_{t,t+H}^{drone, *} &=& \arg\min_{\tilde{\xbf}_{tt}^{drone}, ... , \tilde{\xbf}_{t,t+H}^{drone}} \sum_{t'=t}^{t+H} \tilde{\cbf}_{tt'}^T \Mbf \tilde{\xbf}_{tt'z}^{drone} - (\xbf_{t,t'-1}^{drone})^T\Mbf^T diag(\tilde{\sigma}_{tt'})\Mbf\xbf_{tt'}^{drone}, \\
\label{eqn:linint_constr1}
\text{subject to} &:& \sum_{z\in\Zcal} \tilde{\xbf}_{tt'z}^{drone} = 1  \hspace{4mm} \forall t'\in{t,...,t+H}, \\
\label{eqn:linint_constr2}
& & \Vbf \tilde{\xbf}_{tt'} -  \tilde{\xbf}_{t,t'+1} \geq 0 \hspace{4mm} \forall t'\in{t+1,...,t+H-1}, \\
\label{eqn:linint_constr3}
& & \Vbf \tilde{R}_{tt} -  \tilde{\xbf}_{t,t+1} \geq 0 \hspace{4mm}, \\
\label{eqn:linint_constr4}
& & \tilde{\xbf}_{tt}^{drone}, ... , \tilde{\xbf}_{t,t+H}^{drone} \in \{0,1\}^{|\Zcal|} \hspace{4mm} \forall t'\in{t,...,t+H-1},
\en
where constraint \eqref{eqn:linint_constr1} ensures at each time step the decision chooses a single location, constraints \eqref{eqn:linint_constr2} and \eqref{eqn:linint_constr3} ensure that it can travel to the next location in a single time step, and constraint \eqref{eqn:linint_constr4} keeps the variables binary valued.  This optimization problem has $H|\Zcal|$ variables and can be solved by any software package that can handle quadratic integer programs with linear constraints.  The decision at time $t$ is given by $X^{IE-DLA}(S_t^{drone}) =\tilde{\xbf}_{tt}^{drone, *} $.
\proof{proof.}
See Appendix.
\endproof
\end{lemma}
Using a parameterized interval estimation contribution metric within the lookahead helps the drone choose paths to balance exploration versus exploitation of observations, but also takes into account the requirement to return home which is essential to the success of the mission.

\section{Simulation Results}
\label{sec:results}
In this section, we present the results of the simulation study we conducted in order to determine the best set of joint policies.  We collected real fuel and weather data from an area in Siskiyou county, California to design a realistic wildland fire simulator.  We developed the simulator in python using fuel data gathered from \url{LANDFIRE.org}, and historical wind data gathered from the national oceanic and atmospheric administration (\url{www.ncdc.noaa.gov/}).  The area represents approximately one square mile and is discretized into a set of 30 meter by 30 meter zones.  This case study focuses on fires which have just started and the goal of the agents is to stop the fire from spreading too large before a team of professional firefighters arrive.  Dispatching a set of autonomous agents is less dangerous and quicker than allowing humans to be first responders.  We assume it will take an average of two hours before a team can get to the fire front; hence, the helicopter and drone collaborative team will be simulated until that time.  The decision epochs for both agents are synchronized to occur every 10 minutes because we assume the helicopter must refuel in that time period; hence, $T=12$.  We will measure the performance of the joint policies based on the mean cumulative cost of fire burned throughout the horizon.  We will also consider the probability that the joint policies prevent a class C wildfire (between 10 and 100 acres burned) over the time horizon.  This metric is a heuristic way of computing the risk of the fire growing too large too quickly under a specific set of joint policies.

\begin{figure}[ht]
\includegraphics[width=0.9\textwidth, center]{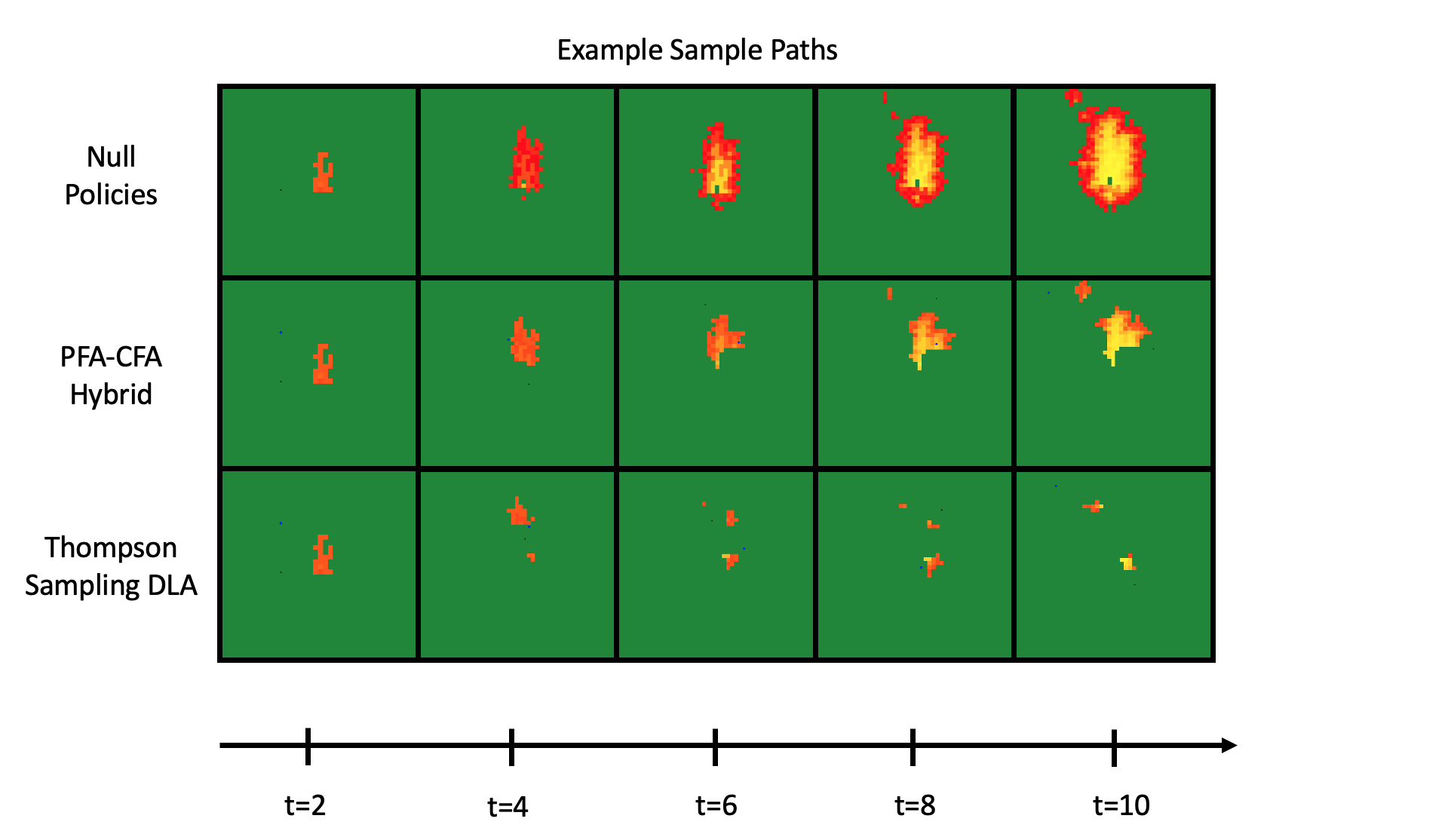}
\caption{Example of a sample path of the environment agent with no intervention (top), a PFA-CFA hybrid learning strategy (middle), and the Thompson sampling direct lookahead learning strategy path (bottom). }
\label{fig:fireseq}
\end{figure}

In the simulator, the agents will begin a sample path at the home base in zone 0, which is the top left corner of the region.  At time $t=0$, there is a randomly generated fire less than 10 acres in area (class B wildfire), and the drone initiates its belief state with a general 25 acre area where each zone has a probability of ignition which decreases radially around the smoke on the horizon.  As the time sequence unfolds, we draw sample paths from the historical wind data to simulate the spread of the fire with Rothermel's fire spread model outlined in \ref{app:Rothermel}.  We also simulate fire spotting by producing a spark ahead of the firefront in the same direction as the wind with some probability.  We assume the drone can accurately make observations 180 meters in each direction, and the helicopter can extinguish a radius of 120 meters with one full load from a given zone.  The helicopter can travel anywhere in the region in one time step, but we limit the drone to traveling 360 meters within one time step.  Figure \ref{fig:fireseq} demonstrates snapshots of a sample path of the environment agent as the fire spreads.

To evaluate which set of policies is the best, any tunable parameters must be optimized with the simulator we've designed.  We've presented the best set of tunable parameters for each set of joint policies in Table \ref{tab:hyper}.

\begin{table}[ht]
{\begin{tabular}{|m{0.75in}||m{1.2in}|m{1.2in}|m{1.2in}|m{1.2in}|}
 \hline
 \multicolumn{5}{|c|}{Optimal Parameter Tuning Results} \\
 \hline
   & Explore & PFA-CFA & IE-DLA ($\theta^{IE})$ & TS-DLA \\
   \hline
   DLA-1 & $(0.1, 2.25) $ & $ (0.1, 2.25)$  & $ (0.1, 2.25, 0.25)$ & $ (0.1, 2.25)$ \\
   \hline
   CFA-DLA ($\theta^{heli}$) & $ (0.1, 2.25, 0) $ & $ (0.2, 2.25, 7)$  & $ (0.1, 2.5, 0.75, 5)$ & $ (0.1, 2.25, 8 )$ \\
   \hline
\end{tabular}}
\caption{Table displaying the optimal parameter tuning results for each set of joint policies.  }
\label{tab:hyper}
\end{table}

\begin{figure}[ht]
\includegraphics[width=.8\textwidth, center]{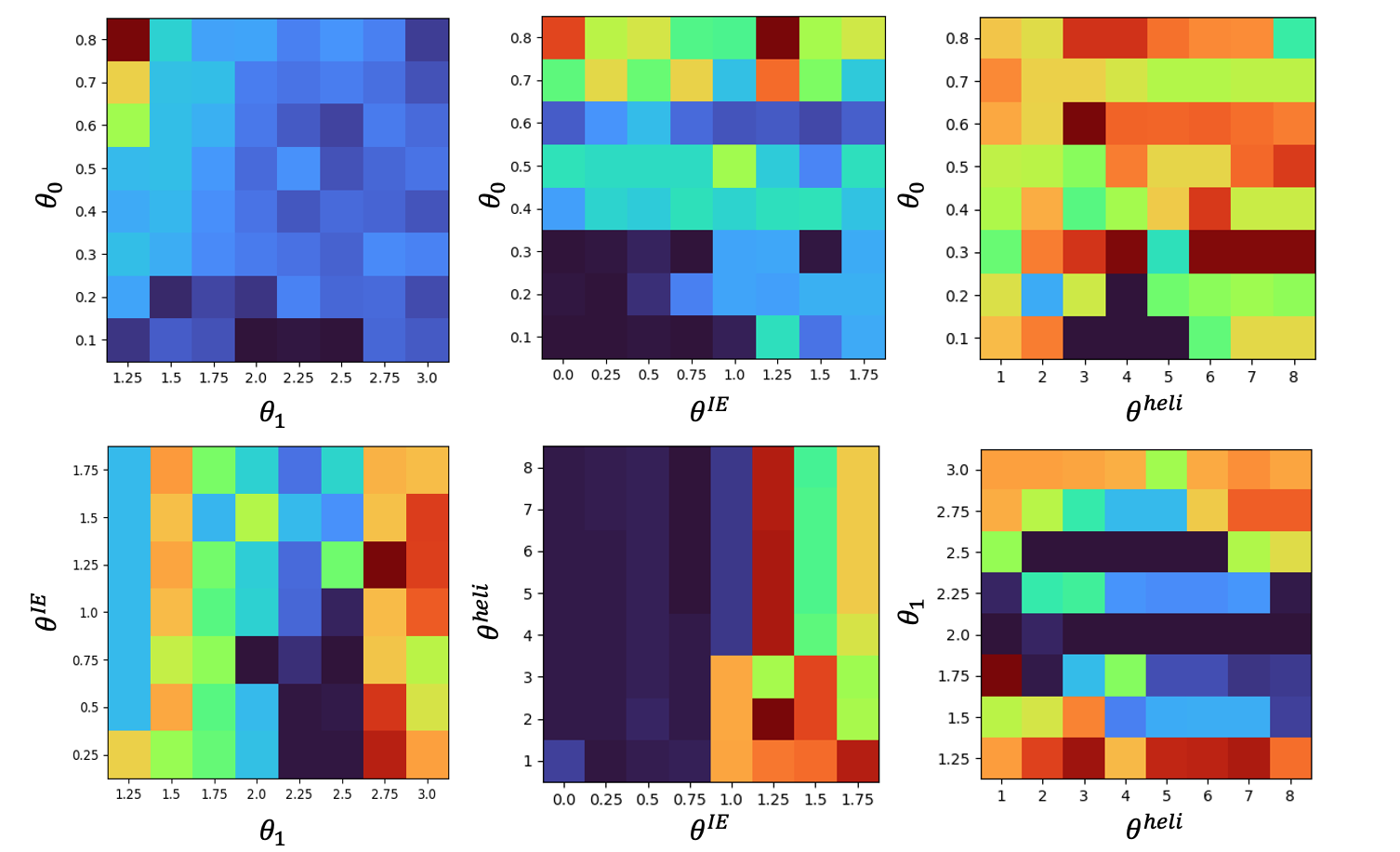}
\caption{The level sets for the interval estimation DLA policy in conjunction with the CFA-DLA helicopter policy. \label{fig:levelsets}}
\end{figure}

We show the level sets of the parameter space near the optimum for the best set of joint policies; the interval estimation DLA in conjunction with the CFA-DLA hybrid.  There are four dimensions to the tunable parameter space and we've plotted a heatmap for each combination of the dimensions in Figure \ref{fig:levelsets}.

\begin{figure}[ht]
\includegraphics[width=1\textwidth, center]{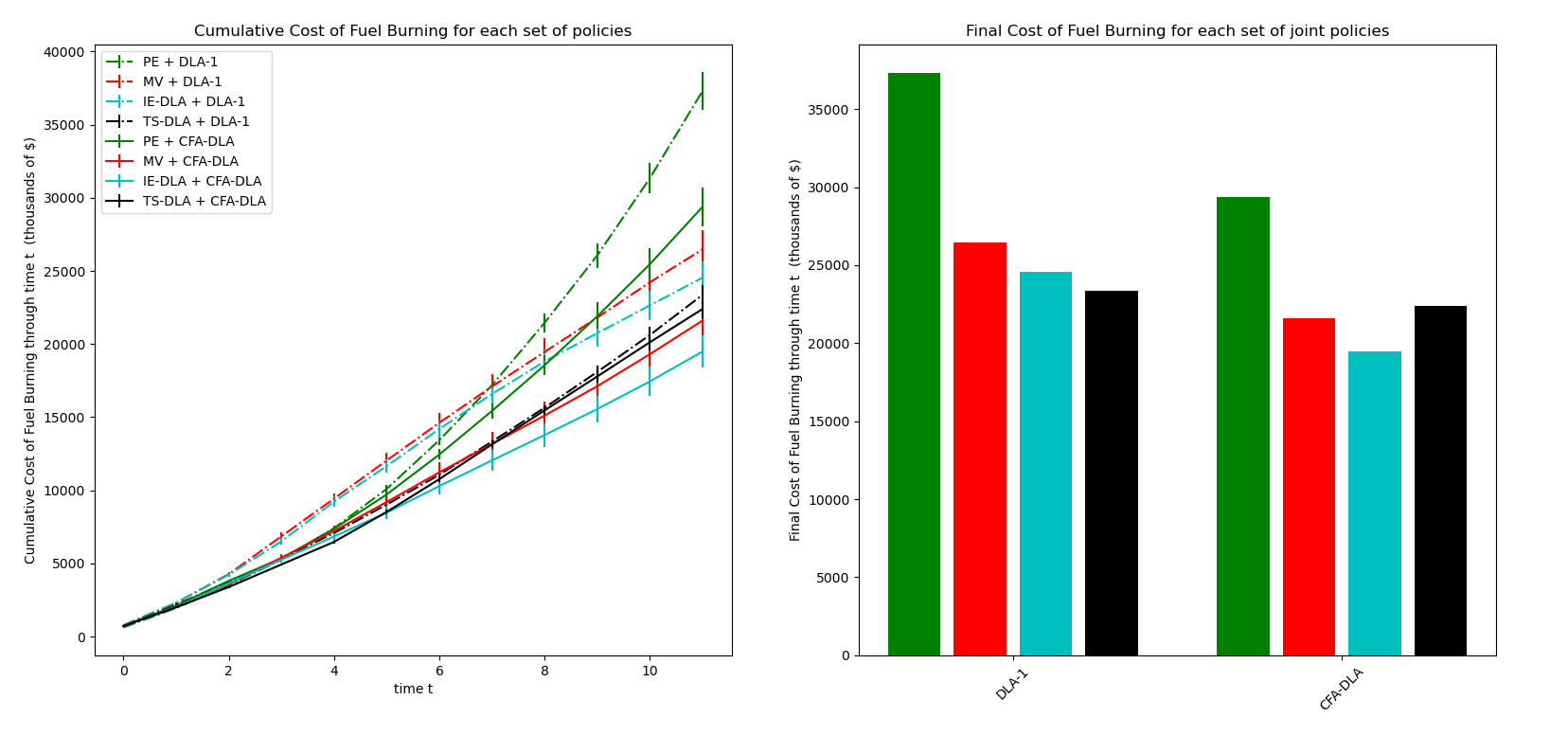}
\caption{The average cumulative cost of burning fuel for each set of joint policies (left) and the final cost of burning fuel for each set of joint policies (right).  }
\label{fig:plots}
\end{figure}

After tuning each of the policies, we evaluated each of the sets of joint policies with the optimal parameter vectors.  In Figure \ref{fig:plots}, we provide the results for the average cumulative cost of burning fuel with confidence intervals for each set of joint policies.  The learning policies in conjunction with the DLA-1 policy are plotted with dotted lines and the CFA-DLA hybrid are plotted with solid lines.

\begin{figure}[ht]
\includegraphics[width=1\textwidth, center]{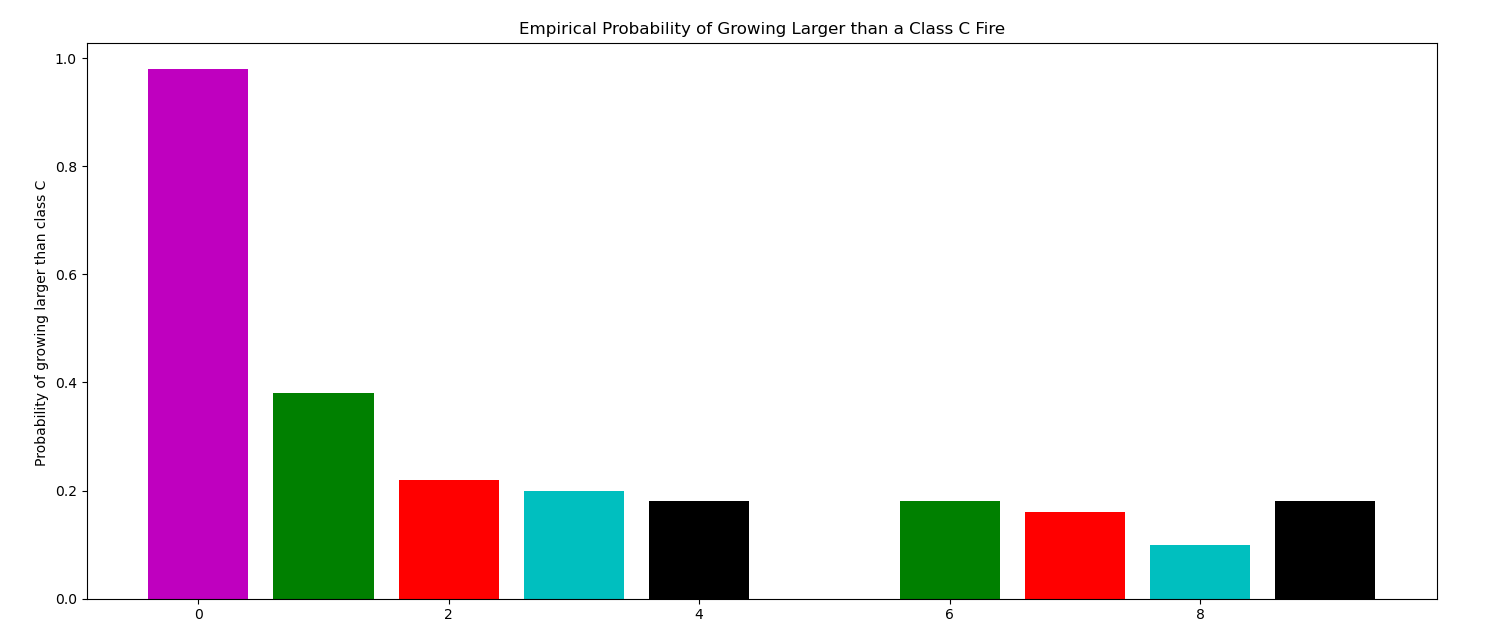}
\caption{The level sets around the optimal tunable parameters in parameter space for the CFA-DLA helicopter policy and the IE-DLA learning policy for the drone. }
\label{fig:pct}
\end{figure}

In Figure \ref{fig:pct}, we provide the probability of the fire growing to be at least a class C wildfire before the time horizon expires. This is a heuristic risk measure we've created to compute the probability of the policy preventing the fire from growing too large before more resources arrive on the scene.  We've included the null policy in this plot to demonstrate that there is a 98\% chance of the fire growing to a class C wildfire when there is no intervention before the fire team arrives.  In this case, the best set of joint policies is the CFA-DLA hybrid and the IE-DLA learning policy which limits the empirical probability of growing to a class C wildfire to 10\%.

From these results, we conclude that the CFA-DLA helicopter policy in conjunction with the Interval Estimation DLA with parameter vector $\theta = (0.1, 2.5, 0.75, 5)$ is the best set of joint policies for minimizing the average cumulative cost of burning fuel and mitigating the risk of large wildfires.

\section{Conclusion}
In this paper, we designed policies for a helicopter equipped with a tanker and a drone to collaboratively extinguish and monitor a wildland fire.  The modeling framework was developed by extending the unified framework for sequential decision problems to multi-agent partially observable systems.  We captured the dynamic uncertainties for each agent with a sophisticated belief model.  Then, we designed policies for both agents which can handle the large state, action, and observation spaces.  The size of the problem presented many computational challenges and limited the set of feasible policies; however, we demonstrate that a class of hybrid parameterized lookahead strategies can be tuned to overcome some of the design challenges to create a robust policy.

After modeling and designing solutions, we created a simulator with real data from a small region in Northern California to discover the best set of joint policies between the two agents.  We discovered the Interval Estimation direct lookahead strategy in conjunction with the hybrid CFA-DLA helicopter policy had the best results in the simulator we designed.

To further extend the work in this paper, one could consider additional physical constraints on the helicopter, the time synchronization assumption between the decision epochs of the agents can be removed, or noise can be added to the communication channel.  Additionally, scenarios with swarms of drones or multiple helicopters can also be considered.  A strength of the multi-agent extension of the unified framework is the feasibility of adding new design considerations into the model.  The standardized organization of the unified framework makes adding new modeling considerations or new agents very straightforward and does not require starting from square one.

\section{Bibliography}


\newpage

\section{Appendix}

\subsection{Approximate Sampling Procedure}
\label{app:sampling}
 This section describes a procedure for drawing samples from the belief model at a time $t$.  At time $t$, we have a belief state given by $B_t = (\pbar_{tz}^K,\pbar_{tz}^H,\pbar_{tz}^Q,\pbar_{tz}^D)_{z\in\Zcal}$.  The estimates of $\pbar_{tz}^K$ for each zone are parameters of a multivariate Bernoulli distribution.  The belief states $\pbar_{tz}^H$, $\pbar_{tz}^Q$, and $\pbar_{tz}^D$ are estimates of the parameters of a multinomial distribution.  We constructed the belief state where a zone is burning if $K_{tz}=1$ and it is completely healthy if $K_{tz} = 0$.  Hence, if we generate samples from the multivariate Bernoulli distribution, then we can generate a sequential sample from the multinomial distribution for the fuel states if the outcome is equal to 1.  If the outcome of the Bernoulli coin flip in each zone is 0, then we assume all of the fuel is healthy.

 The challenge with this sampling procedure is the Bernoulli trials are correlated.  The joint probability distribution for a multivariate Bernoulli distribution has $2^{|\Zcal_{t+1}^{obs}|}-1$ parameters, which would be too computationally intensive for our model.  Hence, we design our sampling procedure using the fitted GPC model from section \ref{sec:belief_model}.  Recall, we assumed the observations were drawn from a function $\sigma(f(z))$, where f(z) is a latent function with a Gaussian process prior and $\sigma(\cdot)$ is a logistic function.  The covariance function for the GP has an RBF structure and the parameterization, $\theta^{cov}$ was already fit by optimizing the log likelihood.  The posterior distribution over the latent function from the previous observations is given by,
\bn
\label{eqn:sample_post}
q(\fbf|\Zcal^{obs}_{t}, (\yhat_z)_{z\in\Zcal^{train}}) = \Ncal\left(\bar{\fbf}, ((\Sigma^{\theta})^{-1} + W)^{-1}\right),
\en
where $W$ is the Hessian of the likelihood function.  However, the vector of means for each zone is no longer relevant because it was blended into our belief model.  Therefore, we assume the vector for the mean of the GP over the latent function is actually the output of the following function of our current belief model,
\bn
\label{eqn:inv_mean}
\bar{\fbf} = \left[\sigma^{-1}\left(\pbar_{tz}^K\right)\right]_{z\in\Zcal}.
\en
If we assume the same correlation structure as the observation model, then we will draw samples through,
\bns
\hat{\fbf} \sim \Ncal\left(\left[\sigma^{-1}\left(\pbar_{tz}^K\right)\right]_{z\in\Zcal}, ((\Sigma^{\theta})^{-1} + W)^{-1}\right) ,
\ens
with equation \eqref{eqn:inv_mean} as the mean of the distribution.  Then we classify the samples for each zone $z\in\Zcal$ through,
\bns
\khat_{tz} = \begin{cases}
    1 & \text{if $\fhat_{z} > 0$},\\
    0 & \text{otherwise}.
\end{cases}
\ens
The samples of the fuel level in each zone are then sampled using the following procedure.  If $\khat_{tz} = 1$, then sample $(\Qhat_{tz}, \Hhat_{tz}, \Dhat_{tz}) \sim Mult(\eta_z, \pbar_{tz}^Q, \pbar_{tz}^H, \pbar_{tz}^D)$.  If $\khat_{tz} = 0$, then $\Hhat_{tz} = \eta_z$ and the remaining cells are zero.  This procedure will provide empirically good samples from the belief model and be computationally efficient.

\subsection{Fire spread simulation model}
\label{app:Rothermel}
The fire spread simulation model was used to construct the adjacency matrix $\hat{\Abf}_{t+1}$ from equation \ref{eqn:env_Kt1_update} for the python simulator.  The adjacency matrix is binary and an entry $\ahat_{t+1,zz'} \in \{0,1\}$ describes whether two zones would connect via ignition if zone $z'$ were ignited at time $t$.  We will draw from the wildfire simulation literature to design a method for computing $\ahat_{t+1,zz'}$ for each entry of the adjacency matrix.  The fuel in the zone determines how far a given flame will spread.  The various characteristics of fuel are given in the following table.

\begin{table}[ht]
{\begin{tabular}{ |p{3.5cm}||p{3.5cm}||p{3.5cm}||p{3.5cm}| }
 \hline
 \multicolumn{4}{|c|}{Fuel Characteristics, $\Xi_z$} \\
 \hline
Variable Name & Fuel Description & Variable Name & Fuel Description\\
\hline
$\kappa_{load}$   & oven dry fuel load & $\kappa_{SAV}$   & surface area to volume ratio \\
\hline
$\kappa_{depth}$   & fuel bed depth & $\kappa_{moist}$   & fuel moisture content \\
\hline
$\kappa_{heat}$   & low heat content & $m_{tot}$   & total mineral content \\
\hline
$\kappa_{dens}$   & oven dry particle density & $m_{eff}$   & effective mineral content \\
\hline
$\lambda_z$   & rate of healthy fuel consumption & $\xi_z$   & rate of natural extinguish \\
\hline
$r_z$   & cost of fuel & $c_z$   & characteristic intensity constant \\
\hline
\end{tabular}}
\caption{Fuel Characteristics within each zone.}
\end{table}

which were drawn from \cite{andrews2018rothermel}.  The set of equations used to measure the rate of spread between two zones is given by the basic fire surface spread model in section 3 of \cite{andrews2018rothermel}, which we will denote $f^{R}(\cdot)$.  Let the set of functions random rate of spread determined by the fuel characteristics for each zone, the wind speed and direction, and the elevation of the land be multiplied by the time step of 10 minutes to compute $\Lhat_{t+1,zz'} = f^{R}(\Xi_z, U_t, \phi_t, \kappa_{elev})$.  If two zones have heterogenous fuel characteristics between them, then we compute the average to generate the random length.  Hence, we compute an entry in the binary matrix, $\hat{a}_{t+1,zz'}$ through,
\bns
\ahat_{t+1,zz'} = \begin{cases}
    1 & \text{if $\Lhat_{t+1,zz'} \geq \|z-z'\|$},\\
    0 & \text{otherwise}.
\end{cases}
\ens
The stochasticity is generated mainly from the wind speed and a phenomenon known as spotting.  Spotting is when an ember jumps out of the fire front (possibly from wind) and ignites regions away from the main flames.  To simulate spotting we sample from a Bernoulli distribution with parameter $p^{spot}$ from every zone where $K_{tz} = 1$. If the Bernoulli random variable is 1, then we generate an exponential random variable with parameter $\lambda^{spot}$ in the wind direction.  This logic will produce a rare tail event which causes a zone to ignite away from the fire.

\subsection{Theoretical Results}
\paragraph{\textbf{Lemma \ref{lem:forecast}}}
\proof{proof.}
Let $\pbar_{tz}^K$ be the estimated probability of a zone being on fire at time $t$ in zone $z$.  Let $\Zcal_{t}^{ext}$ be the set of zones which have been extinguished through time $t$.  Let the length a fire front will travel out of a zone $z'$ in the direction of zone $z$ during the time step between $t$ and $t+1$ be represented by the random variable $\Lhat_{t+1,zz'}$.  The belief state for the zones which are still burning is trivially given by,
\bns
\pbar_{tz}^{K,x} = \pbar_{tz}^K \mathbf{1}_{\Zcal \setminus \Zcal_{t}^{ext}},
\ens
because the zones which have been extinguished are set to zero.  Define the following events as follows:
\bns
E_{tz'} &=& \{\text{The event when zone $z'$ is ignited at time $t$} => K_{tz}^{x} = 1\}, \\
E_{t+1,zz'} &=& \{\text{The event a fire would travel out of zone $z'$ to zone $z$ between $t$ and $t+1$} \}.
\ens
We define $\Pbb[E_{tz'}] = \pbar_{tz}^{K,x}$ and $\Pbb[E_{t+1,zz'}|E_{tz'}] = \Pbb[\Lhat_{t+1,zz'} \geq \|z - z'\| | K_{tz'}=1]$.

The forecasted probability of a zone igniting at time $t+1$ is given by the following measure over the events,
\bns
f^K_{t,t+1,z} = \Pbb\left[\bigcup_{z'\in\Zcal} E_{tz'} \cap E_{t+1,zz'}\right].
\ens
We assume the ignited zones are independent point sources, so if one zone $z'$ would ignite $z$, then it will not impact whether $z''$ would also ignite $z$.  Hence,
\bns
f^K_{t,t+1,z} &=& \Pbb\left[\bigcup_{z'\in\Zcal} E_{tz'} \cap E_{t+1,zz'}\right], \\
&=& 1 - \prod_{z'\in\Zcal} \Pbb\left[(E_{tz'} \cap E_{t+1,zz'})^c\right], \\
&=& 1 - \prod_{z'\in\Zcal}\Pbb\left[ E_{tz'}^c \cup E_{t+1,zz'}^c\right], \\
&=& 1 - \prod_{z'\in\Zcal}\Pbb\left[ E_{tz'}^c\right] + \Pbb\left[E_{t+1,zz'}^c\right] - \Pbb\left[ E_{tz'}^c \cap E_{t+1,zz'}^c\right], \\
&=& 1 - \prod_{z'\in\Zcal} 1 - \pbar_{tz'}^{K,z} + \Pbb\left[E_{t+1,zz'}^c\right] - \Pbb\left[E_{t+1,zz'}^c | E_{tz'}^c \right] \Pbb\left[E_{tz'}^c\right] , \\
&=& 1 - \prod_{z'\in\Zcal} 1 - \pbar_{tz'}^{K,z} + \Pbb\left[E_{t+1,zz'}^c|E_{tz'}\right] \Pbb\left[E_{tz'}\right] + \Pbb\left[E_{t+1,zz'}^c|E_{tz'}^c\right] \Pbb\left[E_{tz'}^c\right]- \Pbb\left[E_{t+1,zz'}^c | E_{tz'}^c \right] \Pbb\left[E_{tz'}^c\right] , \\
&=& 1 - \prod_{z'\in\Zcal} 1 - \pbar_{tz'}^{K,z} +  \Pbb\left[E_{t+1,zz'}^c|E_{tz'}\right] \Pbb\left[E_{tz'}\right] , \\
&=& 1 - \prod_{z'\in\Zcal} 1 - \pbar_{tz'}^{K,z} +  \Pbb[\Lhat_{t+1,zz'} < \|z - z'\| | K_{tz'}=1] \pbar_{tz'}^{K,z},
\ens
where $(\cdot)^c$ represents the set compliment.  This proves equation \eqref{eqn:prior_drift}.
\endproof

\paragraph{\textbf{Lemma \ref{lem:fuel_updates}}}
\proof{proof.}
We break this proof down into a section for observed zones and a section for unobserved zones.  The subsequent paragraph shows the equations for observed zones.

Let $\Ycal_{t+1} = (\yhat_{t+1,z'})_{z\in\Zcal^{obs}_{t+1}}$ be the set of observations from the observed zones.  The model we assumed for the observations was given  by $\yhat_{t+1,z'} = c_z Q_{t+1,z} + \varepsilon_{t+1,z'}$ where $\varepsilon_{t+1,z'} \sim \Ncal(0, \sigma^2)$ is i.i.d. noise with known variance.  If we assume each measurement is independent and the signal to noise ratio is very high, then we can approximate the minimum mean squared error estimate of $Q_{t+1,z}$ through,
\bns
\Qhat_{t+1,z'} &=& \E\left[\frac{Y_{t+1,z'} - \varepsilon_{t+1,z'}}{c_z}|\yhat_{t+1,z'}
\right], \\
&=& \frac{1}{c_z} \E\left[Y_{t+1,z'}|\yhat_{t+1,z'}\right] - \frac{1}{c_z} \E\left[\varepsilon_{t+1,z'}|\yhat_{t+1,z'}\right], \\
&=& \frac{\yhat_{t+1,z'}}{c_z}. \\
\ens
We also make the assumption that any non-burning fuel is healthy because it will always mitigate the risk of fuel that may still be susceptible.

Let a zone $z\in\Zcal \setminus \Zcal_{t+1}^{obs}$ be denoted unobserved at time $t+1$.  The set of unobserved zones have prior belief states from time $t$ which can be forecasted to $t+1$.  They also have some information regarding whether they have ignited using the observation model.  Let $\pbar_{tz}^Q$, $\pbar_{tz}^H$, and $\pbar_{tz}^D$ be the fuel level belief states at time $t$ in zone $z$.  Let $\pbar_{t+1,z}^K$ represent the belief about $K_{t+1,z}$ we computed in equation \eqref{eqn:GPC_pred}.  Before we forecast the fuel levels, we must compute the forecast of whether fire will start in a zone or if it is still burning.  We denote $\E_t[\cdot|B_t, \xbf_{t}^{heli}] = \E_t[\cdot]$.   Respectively, the forecasted event while at time $t$ that a fire will start in zone $z$ at time $t+1$ and the event that a fire will continue to burn in zone $z$ at time $t+1$, are given by,
\bns
E^{start}_{t+1,z} &=& \{K_{t+1,z} = 1\} \cap \{K_{tz}^x = 0\},\\
E^{burn}_{t+1,z} &=& \{K_{t+1,z} = 1\} \cap \{K_{tz}^x = 1\}.\\
\ens
Based on the assumptions from section \ref{sec:problem_desc}, we define the following conditional probabilities,
\bns
\Pbb[K_{tz}^{x} = 0|K_{t+1,z} = 0] &=& 1,  \\
\Pbb[K_{tz}^{x} = 1|K_{t+1,z} = 0] &=& 0.
\ens
Hence, the forecasted probability measure over the event the fire will ignite at $t+1$ is given the belief states is given by,
\bns
f_{t,t+1,z}^{start} &=& \E_t[E^{start}_{t+1,z}], \\
 &=& \E_t[\{K_{t+1,z} = 1\} \cap \{K_{tz}^x = 0\}], \\
 &=& \Pbb\left[K_{tz}^x = 0|K_{t+1,z} = 1\right] \Pbb\left[K_{t+1,z} = 1\right], \\
&=& \Pbb\left[K_{tz}^x = 0\right] - \Pbb\left[K_{tz}^x = 0|K_{t+1,z} = 0\right] \Pbb\left[K_{t+1,z} = 0\right], \\
&=& \Pbb\left[K_{tz}^x = 0\right] - \Pbb\left[K_{t+1,z} = 0\right], \\
&=& (1 - \pbar_{tz}^{K,x}) - (1-f_{t,t+1,z}^K) =  f_{t,t+1,z}^K - \pbar_{tz}^{K,x}.\\
\ens
The forecasted probability measure over the event the fire is continuing to burn at $t+1$ is given the belief states is given by,
\bns
f_{t,t+1,z}^{burn} &=& \E_t[E^{burn}_{t+1,z}], \\
&=& \E_t[\{K_{t+1,z} = 1\} \cap \{K_{tz}^x = 1\}], \\
 &=& \Pbb\left[K_{tz}^x = 1|K_{t+1,z} = 1\right] \Pbb\left[K_{t+1,z} = 1\right], \\
&=& \Pbb\left[K_{t+1,z} = 1\right] - \Pbb\left[K_{tz}^x = 0|K_{t+1,z} = 1\right] \Pbb\left[K_{t+1,z} = 1\right], \\
&=& \pbar_{t+1,z}^K - \pbar_{t+1,z}^K + \pbar_{tz}^{K,x}, \\
&=& \pbar_{tz}^{K,x}. \\
\ens

After new observations, the Bayesian update is given by,
\bns
\Pbb[E_{t+1,z}^{start} | \Ycal_{t+1}, B_t, \xbf_{tz}^{heli}] &=& \frac{\Pbb[\Ycal_{t+1} | E^{start}_{t+1,z}, B_t, \xbf_{tz}^{heli}] \Pbb[E_{t+1,z}^{start} | B_t, \xbf_{tz}^{heli}]}{\Pbb[\Ycal_{t+1} | B_t, \xbf_{tz}^{heli}]}, \\
\Pbb[E_{t+1,z}^{burn} | \Ycal_{t+1}, B_t, \xbf_{tz}^{heli}] &=& \frac{\Pbb[\Ycal_{t+1} | E^{burn}_{t+1,z}, B_t, \xbf_{tz}^{heli}] \Pbb[E_{t+1,z}^{burn} | B_t, \xbf_{tz}^{heli}]}{\Pbb[\Ycal_{t+1} | B_t, \xbf_{tz}^{heli}]}.
\ens
However, the new observations do not provide any new information about the events in question.  Hence the following equality holds,
\bns
\frac{\Pbb[\Ycal_{t+1} | E^{start}_{t+1,z}, B_t, \xbf_{tz}^{heli}]}{\Pbb[\Ycal_{t+1} | B_t, \xbf_{tz}^{heli}]} = \frac{\Pbb[\Ycal_{t+1} | E^{burn}_{t+1,z}, B_t, \xbf_{tz}^{heli}]}{\Pbb[\Ycal_{t+1} | B_t, \xbf_{tz}^{heli}]}.
\ens
Since the events are disjoint, the previous equality also implies the ratio of each of the posterior to the prior will be the same for each event.  Hence
\bns
\frac{\Pbb[E_{t+1,z}^{start} \bigcup E_{t+1,z}^{burn} | \Ycal_{t+1}, B_t, \xbf_{tz}^{heli}]}{\Pbb[E_{t+1,z}^{start} \bigcup E_{t+1,z}^{burn} | B_t, \xbf_{tz}^{heli}]} = \frac{\Pbb[E_{t+1,z}^{start} | \Ycal_{t+1}, B_t, \xbf_{tz}^{heli}]}{\Pbb[E_{t+1,z}^{start} | B_t, \xbf_{tz}^{heli}]} = \frac{\Pbb[E_{t+1,z}^{burn} | \Ycal_{t+1}, B_t, \xbf_{tz}^{heli}]}{\Pbb[E_{t+1,z}^{burn} | B_t, \xbf_{tz}^{heli}]},
\ens
where $\Pbb[E_{t+1,z}^{start} \bigcup E_{t+1,z}^{burn} | \Ycal_{t+1}, B_t, \xbf_{tz}^{heli}] = \pbar_{t+1,z}^K$, $\Pbb[E_{t+1,z}^{start} \bigcup E_{t+1,z}^{burn} | B_t, \xbf_{tz}^{heli}] = f_{t,t+1,z}^K$,  $\Pbb[E_{t+1,z}^{start} | B_t, \xbf_{tz}^{heli}] = f_{t,t+1,z}^K - \pbar_{tz}^{K,x}$, and $\Pbb[E_{t+1,z}^{burn} | B_t, \xbf_{tz}^{heli}] = \pbar_{tz}^{K,x}$.  After rearranging the terms from the equations above, the following posteriors, $\pbar_{t+1,z}^{start}$ and $\pbar_{t+1,z}^{burn}$, are given by,
\bns
\pbar_{t+1,z}^{start} &=& \pbar_{t+1,z}^K -  \frac{\pbar_{t+1,z}^K \pbar_{tz}^{K,x}}{f_{t,t+1,z}},\\
\pbar_{t+1,z}^{start} &=& \frac{\pbar_{t+1,z}^K \pbar_{tz}^{K,x}}{f_{t,t+1,z}}.
\ens
Since these represent probabilities, and $\pbar_{tz}^{K,x} \leq f_{t,t+1,z}^K$, then the terms will always be between 0 and 1 unless $f_{t,t+1,z} = 0$ and $\pbar_{tz}^{K,x} = 0$.  In the special case where the terms become undefined, then we assert $\pbar_{t+1,z}^{start} = \pbar_{t+1,z}^K$ and $\pbar_{t+1,z}^{burn} = 0$.

We compute the fuel level forecast by taking the conditional expectation of equations \eqref{eqn:env_Qtx_update} - \eqref{eqn:env_Dt1_update} with respect to the belief state at time $t$ and observations at $t+1$. The conditional expectation of equation \eqref{eqn:env_Qt1_update} for all unobserved zones is given by,
\bns
\pbar_{t+1,z}^Q &=&  \E_t\left[\frac{Q_{t+1}}{\eta_z}\right],\\
&=&  \frac{1}{\eta_z}\E_t\left[Q_{tz}^x - \xi_z Q_{tz}^x + \lambda_z Q_{tz}^x H_{tz} + Q^{init}_{z} (K_{t+1,z} - K^x_{tz})\right],
\ens
which is equivalent to,
\bns
\pbar_{t+1,z}^Q &=&  \frac{1}{\eta_z} \left(\E_t\left[Q_{tz}E^{burn}_{t+1,z}\right] - \E_t\left[\xi_z Q_{tz}E^{burn}_{t+1,z}\right] + \E_t\left[\lambda_z Q_{tz} H_{tz}E^{burn}_{t+1,z}\right] + \E_t\left[Q^{init}_{z} E^{start}_{t+1,z}\right]\right).\\
\ens
Given the belief state at time $t$, the fuel level is independent from the events determining whether the zone just started burning or not.  Hence, these expectations can be separated, given by,
\bns
\pbar_{t+1,z}^Q &=&  \frac{1}{\eta_z} \left( \E_t\left[Q_{tz}\right]\E_t\left[E^{burn}_{t+1,z}\right] - \xi_z\E_t\left[ Q_{tz}\right]\E_t\left[E^{burn}_{t+1,z}\right] + \lambda_z \E_t\left[ Q_{tz} H_{tz}\right]\E_t\left[E^{burn}_{t+1,z}\right] + Q^{init}_{z}\E_t\left[ E^{start}_{t+1,z}\right]\right),\\
&=&  \frac{1}{\eta_z} \left( \E_t\left[Q_{tz}\right] \pbar_{tz}^{K,x} - \xi_z\E_t\left[ Q_{tz}\right] \pbar_{tz}^{K,x} + \lambda_z \E_t\left[ Q_{tz} H_{tz}\right] \pbar_{tz}^{K,x} + Q^{init}_{z} \pbar_{t+1,z}^{start}\right).
\ens
Recall, the belief state variables for the fuel levels are assumed to be drawn from a multinomial distribution.  The following moments of the multinomial distributions are substituted into the expression,
\bns
\E_t[Q_{tz}] &=& \eta_z \pbar_{tz}^Q, \\
\E_t[H_{tz}] &=& \eta_z \pbar_{tz}^H, \\
\E_t[D_{tz}] &=& \eta_z \pbar_{tz}^D, \\
\E_t[Q_{tz}H_{tz}] &=& \eta_z (\eta_z - 1) \pbar_{tz}^Q\pbar_{tz}^H.
\ens
The updated belief state for equation \eqref{eqn:env_Qt1_update} is given by,
\bns
\pbar_{t+1,z}^Q &=&  \frac{1}{\eta_z} \left( \eta_z \pbar_{tz}^Q \pbar_{tz}^{K,x} - \xi_z\eta_z \pbar_{tz}^Q \pbar_{tz}^{K,x} + \lambda_z \eta_z (\eta_z - 1) \pbar_{tz}^Q\pbar_{tz}^H \pbar_{tz}^{K,x} + Q^{init}_{z} \pbar_{t+1,z}^{start}\right),\\
&=& \pbar_{tz}^Q \pbar_{tz}^{K,x} - \xi_z \pbar_{tz}^Q \pbar_{tz}^{K,x} + \lambda_z (\eta_z - 1) \pbar_{tz}^Q\pbar_{tz}^H \pbar_{tz}^{K,x} + \frac{Q^{init}_{z}}{\eta_z} \pbar_{t+1,z}^{start},\\
&=& \left(1 - \xi_z + \lambda_z (\eta_z - 1) \pbar_{tz}^H \right) \pbar_{tz}^Q \pbar_{tz}^{K,x} + \frac{Q^{init}_{z}}{\eta_z} \pbar_{t+1,z}^{start}.
\ens
A similar procedure is performed for the conditional expectations of equation \eqref{eqn:env_Ht1_update} and \eqref{eqn:env_Dt1_update} .  We outline the steps for each of the derivations in the subsequent sets of equations, given by,
\bns
\pbar_{t+1,z}^H &=&  \E_t\left[\frac{H_{t+1}}{\eta_z}\right],\\
&=& \frac{1}{\eta_z}\E_t\left[H_{tz} - \lambda_z Q_{tz} H_{tz}E^{burn}_{t+1,z} - Q^{init}_{z} E^{start}_{t+1,z}\right], \\
&=& \frac{1}{\eta_z} \left(\E_t\left[H_{tz}\right] - \lambda_z \E_t\left[Q_{tz} H_{tz}\right]\E_t\left[E^{burn}_{t+1,z}\right] - Q^{init}_{z} \E_t\left[E^{start}_{t+1,z}\right] \right), \\
&=& \pbar_{tz}^H - \lambda_z (\eta_z-1)\pbar_{tz}^Q\pbar_{tz}^H \pbar_{tz}^{K,x} - Q^{init}_{z} \pbar_{t+1,z}^{start}, \\
\pbar_{t+1,z}^D &=&  \E_t\left[\frac{D_{t+1}}{\eta_z}\right],\\
&=& \frac{1}{\eta_z}\E_t\left[ D_{tz} + Q_{tz} (1 - K_{tz}^x) + \xi_z Q^x_{tz} \right], \\
&=& \frac{1}{\eta_z} \left( \E_t\left[ D_{tz} \right] + \E_t\left[Q_{tz}\right] - \E_t\left[ Q_{tz} E^{burn}_{t+1,z}\right] + \xi_z \E_t\left[ Q_{tz} E^{burn}_{t+1,z}\right] \right), \\
&=& \pbar_{tz}^D + \pbar_{tz}^{Q} (1 - \pbar_{tz}^{K,x}) + \xi_z \pbar_{tz}^{Q} \pbar_{tz}^{K,x}. \\
\ens

\endproof

\paragraph{\textbf{Lemma \ref{lem:samp_LA}}}
\proof{proof.}
Let $\tilde{\xbf}_{tt}^{drone}$, $\tilde{\xbf}_{t,t+1}^{drone}$, ... , $\tilde{\xbf}_{t,t+H}^{drone}$ be the binary vectors representing the drone decisions $H$ steps into the lookahead model.  Let $\Mbf \in \{0,1\}^{|\Zcal| \times |\Zcal|}$ be a matrix where each column is a binary mask over the field of view for each possible location in the set of zones.  Let $\Vbf \in \{0,1\}^{|\Zcal| \times |\Zcal|}$ be a matrix where each column is a binary mask over the zones reachable within 10 minutes for each possible location in the set of zones.  The objective function is given by,
\bns
\max_{\tilde{\xbf}_{tt}^{drone}\in \tilde{\Xcal}_{tt}^{drone}, ... , \tilde{\xbf}_{t,t+H}^{drone}\in\tilde{\Xcal}_{t,t+H}^{drone}} \frac{1}{m} \sum_{i=1}^{m} \tilde{C}\left(\Stilde_{tt}^{drone}(\tilde{\omega}^{i}), \tilde{\xbf}_{tt}^{drone}\right) +  \sum_{t'=t+1}^{t+H} \tilde{C}\left(\Stilde_{tt'}^{drone}(\tilde{\omega}^{i}), \tilde{\xbf}_{tt'}^{drone}\right),
\ens
where, $\tilde{C}\left(\Stilde_{tt'}^{drone}, \tilde{\xbf}_{tt'}^{drone}\right) = \sum_{z\in\Zcal_{tt'}^{obs}( \tilde{\xbf}_{tt'z}^{drone})} r_z \left(\tilde{Q}_{t,t'+1,z} - \tilde{Q}_{ttz}^x\right) - c^{fail} \mathbb{I}_{tt'}$.  The binary vector representing the mask over the observed zones for a given decision is given by $\Mbf \tilde{\xbf}_{tt'}^{drone}$.  Hence, if we define $\tilde{\cbf}_{tt'} = \left[\frac{1}{m} \sum_{i=1}^{m} r_z (\tilde{Q}_{t,t'+1,z}(\tilde{\omega}^{i}) - \tilde{Q}_{tt'z}^x(\tilde{\omega}^{i}))\right]_{z\in\Zcal} - c^{fail} \mathbb{I}_{tt'}$ to be the vector of costs.  Then, the single period contribution function decomposes to,
\bns
\tilde{C}\left(\Stilde_{tt'}^{drone}, \tilde{\xbf}_{tt'}^{drone}\right) = \tilde{\cbf}_{tt'}^T \Mbf \tilde{\xbf}_{tt'}^{drone}.
\ens
After rearranging the objective function with the vectorization presented, we arrive at,
\bns
&\max_{\tilde{\xbf}_{tt}^{drone}\in \tilde{\Xcal}_{tt}^{drone}, ... , \tilde{\xbf}_{t,t+H}^{drone}\in\tilde{\Xcal}_{t,t+H}^{drone}}& \tilde{\cbf}_{tt}^T \Mbf \tilde{\xbf}_{tt}^{drone} + \sum_{t'=t+1}^{t+H} \tilde{\cbf}_{tt'}^T \Mbf \tilde{\xbf}_{tt'}^{drone},\\
&\max_{\tilde{\xbf}_{tt}^{drone}\in \tilde{\Xcal}_{tt}^{drone}, ... , \tilde{\xbf}_{t,t+H}^{drone}\in\tilde{\Xcal}_{t,t+H}^{drone}}& \sum_{t'=t}^{t+H} \tilde{\cbf}_{tt'}^T \Mbf \tilde{\xbf}_{tt'}^{drone}.
\ens

The decision space for the decision at time $t$ is given by,
\bns
\tilde{\Xcal}_{tt'}^{drone} = \left\{ \tilde{\xbf}_{tt'}^{drone} \in \{0,1\}^{|\Zcal|} : \sum_{z\in\Zcal} \xtilde_{tt'z}^{drone} = 1,  \Vbf \tilde{\xbf}_{tt'}^{drone} - \tilde{R}_{tt'}^{drone} \geq \mathbf{0}\right\},
\ens
where $\mathbf{0}$ is the $|\Zcal|$-dimensional vector of zeros.  The first constraint allows only one zone to be chosen and the second constraint allows the drone to travel within its maximum distance for one step.  Hence, our objective function becomes equation \eqref{eqn:samp_LA} and the constraints are given by \eqref{eqn:linint_constr1s}, \eqref{eqn:linint_constr2s}, \eqref{eqn:linint_constr3s}. There are $H$ decision vectors with $|\Zcal|$ variables; hence, the optimization problem has $H|\Zcal|$ binary variables.  The objective function is a sum of inner products between variables; hence, it is linear.  \eqref{eqn:linint_constr1s} and \eqref{eqn:linint_constr2s} are linear functions.  \eqref{eqn:linint_constr3s} restricts the variables to be binary valued.  Hence, the problem reduces to a  integer linear program which can be solved by any software package with integer linear values.  In this problem, we are only interested in the next decision, so the policy chooses the solution to the problem ,$\xbf_{tt}^{drone,*}$.
\endproof

\paragraph{\textbf{Lemma \ref{lem:prob_LA}}}
\proof{proof.}
Let $\tilde{\xbf}_{tt}^{drone}$, $\tilde{\xbf}_{t,t+1}^{drone}$, ... , $\tilde{\xbf}_{t,t+H}^{drone}$ be the binary vectors representing the drone decisions $H$ steps into the lookahead model.  Let $\Mbf \in \{0,1\}^{|\Zcal| \times |\Zcal|}$ be a matrix where each column is a binary mask over the field of view for each possible location in the set of zones.  Let $\Vbf \in \{0,1\}^{|\Zcal| \times |\Zcal|}$ be a matrix where each column is a binary mask over the zones reachable within 10 minutes for each possible location in the set of zones.  The objective function is given by,
\bns
\max_{\tilde{\xbf}_{tt}^{drone}\in \tilde{\Xcal}_{tt}^{drone}, ... , \tilde{\xbf}_{t,t+H}^{drone}\in\tilde{\Xcal}_{t,t+H}^{drone}}  \tilde{C}\left(\Stilde_{tt}^{drone}, \tilde{\xbf}_{tt}^{drone}\right) +  \sum_{t'=t+1}^{t+H} \tilde{C}\left(\Stilde_{tt'}^{drone}, \tilde{\xbf}_{tt'}^{drone}\right),
\ens
where, $\tilde{C}\left(\Stilde_{tt'}^{drone}, \tilde{\xbf}_{tt'}^{drone}\right) = \sum_{z\in\Zcal_{tt'}^{obs}( \tilde{\xbf}_{tt'z}^{drone})} \tilde{\mu}_{tt'z} + \theta^{IE} \tilde{\sigma}_{tt'z}$.  The binary vector representing the mask over the observed zones for a given decision is given by $\Mbf \tilde{\xbf}_{tt'}^{drone}$.  Then, the single period contribution function decomposes to,
\bns
\tilde{C}\left(\Stilde_{tt'}^{drone}, \tilde{\xbf}_{tt'}^{drone}\right) &=& \left(\tilde{\mu}_{tt'} + \theta^{IE} \tilde{\sigma}_{tt'} \odot \left(\left(\mathbf{1}\mathbf{1}^T - \Mbf\right) \Rtilde_{tt'}^{drone}\right)\right)^T   \Mbf \tilde{\xbf}_{tt'}^{drone},\\
&=& \left(\tilde{\mu}_{tt'} + \theta^{IE} diag(\tilde{\sigma}_{tt'}) \mathbf{1}\mathbf{1}^T\Rtilde_{tt'}^{drone}  -\theta^{IE} diag(\tilde{\sigma}_{tt'}) \Mbf \Rtilde_{tt'}^{drone}\right)^T   \Mbf \tilde{\xbf}_{tt'}^{drone}, \\
&=& \left(\tilde{\mu}_{tt'} + \theta^{IE} \tilde{\sigma}_{tt'}\right)^T \Mbf \tilde{\xbf}_{tt'}^{drone}  -\theta^{IE} (\Mbf\Rtilde_{tt'}^{drone})^T diag(\tilde{\sigma}_{tt'})\Mbf \tilde{\xbf}_{tt'}^{drone}.
\ens
After rearranging the objective function with the vectorization presented above, we arrive at,
\bns
&\max_{\tilde{\xbf}_{tt}^{drone}\in \tilde{\Xcal}_{tt}^{drone}, ... , \tilde{\xbf}_{t,t+H}^{drone}\in\tilde{\Xcal}_{t,t+H}^{drone}}& \sum_{t'=t}^{t+H} \left(\tilde{\mu}_{tt'} + \theta^{IE} \tilde{\sigma}_{tt'}\right)^T \Mbf \tilde{\xbf}_{tt'}^{drone}  -\theta^{IE} (\Mbf\Rtilde_{tt'}^{drone})^T diag(\tilde{\sigma}_{tt'})\Mbf \tilde{\xbf}_{tt'}^{drone}.
\ens
Recall, $\Rtilde_{tt'}^{drone} = \tilde{\xbf}_{t,t'-1}^{drone}$, and let $\tilde{\cbf}_{tt'} = \tilde{\mu}_{tt'} + \theta^{IE} \tilde{\sigma}_{tt'}$; hence, we can write the objective function as a quadratic function given by,
\bns
\max_{\tilde{\xbf}_{tt}^{drone}\in \tilde{\Xcal}_{tt}^{drone}, ... , \tilde{\xbf}_{t,t+H}^{drone}\in\tilde{\Xcal}_{t,t+H}^{drone}}\tilde{\cbf}_{tt'}^T \Mbf \tilde{\xbf}_{tt'}^{drone} & -& \theta^{IE} (\Mbf\tilde{\Rtilde}_{tt}^{drone})^T diag(\tilde{\sigma}_{tt'})\Mbf \tilde{\xbf}_{tt'}^{drone} \\ &+& \sum_{t'=t+1}^{t+H} \tilde{\cbf}_{tt'}^T \Mbf \tilde{\xbf}_{tt'}^{drone}  - \theta^{IE} (\Mbf\tilde{\xbf}_{t,t'-1}^{drone})^T diag(\tilde{\sigma}_{tt'})\Mbf \tilde{\xbf}_{tt'}^{drone}.
\ens

The decision space for the decision at time $t$ is given by,
\bns
\tilde{\Xcal}_{tt'}^{drone} = \left\{ \tilde{\xbf}_{tt'}^{drone} \in \{0,1\}^{|\Zcal|} : \sum_{z\in\Zcal} \xtilde_{tt'z}^{drone} = 1,  \Vbf \tilde{\xbf}_{tt'}^{drone} - \tilde{R}_{tt'}^{drone} \geq \mathbf{0}\right\},
\ens
where $\mathbf{0}$ is the $|\Zcal|$-dimensional vector of zeros.  The first constraint allows only one zone to be chosen and the second constraint allows the drone to travel within its maximum distance for one step.  Hence, our objective function becomes equation \eqref{eqn:samp_LA} and the constraints are given by \eqref{eqn:linint_constr1s}, \eqref{eqn:linint_constr2s}, \eqref{eqn:linint_constr3s}. There are $H$ decision vectors with $|\Zcal|$ variables; hence, the optimization problem has $H|\Zcal|$ binary variables.  The objective function is a sum of inner products between variables; hence, it is linear.  \eqref{eqn:linint_constr1s} and \eqref{eqn:linint_constr2s} are linear functions.  \eqref{eqn:linint_constr3s} restricts the variables to be binary valued.  Hence, the problem reduces to a  integer linear program which can be solved by any software package with integer linear values.  In this problem, we are only interested in the next decision, so the policy chooses the solution to the problem ,$\xbf_{tt}^{drone,*}$.
\endproof

\end{document}